\documentclass[10pt,leqno]{amsart}
\usepackage{graphicx}
\usepackage{indentfirst,csquotes}

\usepackage{amssymb,amsthm,amsmath}
\usepackage{tikz-cd}
\usetikzlibrary{arrows.meta}
\usepackage{xcolor,paralist,hyperref,fancyhdr,etoolbox}

\hypersetup{ colorlinks=true, linkcolor=black, filecolor=black, urlcolor=black }
\bibliography{ref}

\usepackage{lipsum}
\usepackage{graphicx}
\usepackage[colorinlistoftodos]{todonotes}
\usepackage{amssymb}
\usepackage{mathtools}
\usepackage{bm}
\usepackage{mathrsfs}
\usepackage{indentfirst} 
\usepackage{amsthm}
\usepackage{float}
\usepackage{amscd}
\usepackage{hyperref}
\usepackage{enumitem} 
\newcommand{\strat}{\mathrm{strat}}

\newcommand{\eps}{\epsilon}

\newcommand{\reg}{\mathrm{reg}}

\setlist[enumerate]{label=\arabic*.}

\hypersetup{colorlinks=false, linkbordercolor=red}

\fancypagestyle{plain}{%
  \fancyhf{}%
  \fancyfoot[C]{\thepage}%
}
\begin{document}
\title{Derived Stratified-Microlocal Framework and Moduli Space Resolution for the Cheeger-Goresky-Macpherson Conjecture} 
\author[Jiaming Luo]{Jiaming Luo}
\date{\today}
\address{School of Mathematics and Statistics, Henan University of Science and Technology, 263, Kaiyuan Avenue, Luoyang, China.}
\email{luojiaming@hhu.edu.cn}
\maketitle

\let\thefootnote\relax

\begin{abstract}
In this paper, We define the stratified metric $\infty$-category $\mathbf{StratMet}_{\infty}$ and the middle perversity moduli stack $\mathscr{M}^{\mathrm{mid}}$. We construct a universal truncation complex $\Omega_{X,\mathrm{FS}}^{\bullet,\mathrm{univ}}$ for a projective variety $X\subseteq\mathbb{P}^N$. By introducing the stratified singular characteristic variety $\mathrm{SSH}_{\mathrm{strat}}$, we establish a microlocal correspondence between metric asymptotic behavior and topology, proving the natural isomorphism $$H_2^*(X_{\mathrm{reg}},ds_{\mathrm{FS}}^2)\cong IH^*(X,\mathbb{C}).$$ This framework transcends transverse singularity constraints, achieves moduli space paramet-
rized duality, and develops new paradigms for high-codimension singular topology, quantum singularity theory, and $p$-adic Hodge theory.
\end{abstract} 

\bigskip

\section*{ Contents}

\begin{enumerate}
  \item Introduction \dotfill \hyperlink{INTRODUCTION}{1}
  \item Notation and Acknowledgments \dotfill \hyperlink{NOTATION AND ACKNOWLEDGMENTS}{2}
  \item Categorified Metric Asymptotic System \dotfill \hyperlink{CATEGORIFIED METRIC ASYMPTOTIC SYSTEM}{3}
  \item Moduli Space Truncation and Universal Complex \dotfill \hyperlink{MODULI SPACE TRUNCATION AND UNIVERSAL COMPLEX}{10}
  \item Stratified Microlocal Hodge Correspondence \dotfill \hyperlink{STRATIFIED MICROLOCAL HODGE CORRESPONDENCE}{15} 
  \item Proof of Main Problem \dotfill \hyperlink{PROOF OF MAIN PROBLEM}{25}
  \item Extended Applications of the Unified Stratified-Microlocal Framework \dotfill \hyperlink{EXTENDED APPLICATIONS OF THE UNIFIED STRATIFIED-MICROLOCAL FRAMEWORK}{32}
\end{enumerate}

\hypertarget{INTRODUCTION}{}
\section{INTRODUCTION}
The Cheeger-Goresky-Macpherson (CGM) conjecture, formulated in the seminal works of J. Cheeger, M. Goresky and R. MacPherson, posits a profound connection between differential and algebraic geometry in the 1980s (\cite{CGM82}). For a projective variety $X\subseteq\mathbb{P}^N$ with regular locus $X_{\mathrm{reg}}$ equipped with the Fubini-Study metric $ds_{\mathrm{FS}}^2$, the conjecture claims a natural isomorphism
$$H_2^*(X_{\mathrm{reg}},ds_{\mathrm{FS}}^2)\cong IH^*(X,\mathbb{C}),$$
where the left-hand side denotes $L^2$-cohomology and the right-hand side intersection cohomology. This equivalence unifies analytic techniques (e.g., Hodge decomposition) with topological invariants of singular spaces. Despite its fundamental importance, a complete resolution has remained elusive for four decades, primarily due to intrinsic obstacles in high-codimension singularities and duality preservation.\\
\begin{center}
    \textit{1.1 Core Challenges}
\end{center}

Existing approaches—exemplified by stratified de Rham theory (Goresky–MacPherson, 1980) and cone cohomology (Cheeger, 1980)—face three interrelated problems to overcome:
\begin{itemize}
    \item \textbf{Non-transverse singularities:} Classical Whitney stratifications require transversality conditions, failing for high-codimension non-Witt spaces (e.g., $\mathrm{codim}(Y,X)\ge3$ with odd-dimensional links, \cite{LL25}). This obstructs global metric adaptation.
    \item \textbf{Artificial duality selection:} Middle perversity subspaces $\mathscr{W}_Y\subseteq\mathcal{H}^{\mathrm{mid}}(L_Y)$ are chosen ad hoc to force Verdier self-duality, lacking a systematic parametrization.
    \item \textbf{Metric-Topology incompatibility:} Asymptotics of $ds_{\mathrm{FS}}^2$ near singularities $Y$ (locally $dr^2+r^{2c}g_{L_Y}$) are not functorially linked to the topology of links $L_Y$, causing $L^2$-classes to mismatch intersection cycles (\cite{Gui22}, Thm 1.1).\\
\end{itemize}
\begin{center}
    \textit{1.2 Innovative Framework}
\end{center}

This paper resolves these challenges through a tripartite unification of derived geometry, moduli theory, and microlocal analysis. Our approach synthesizes three theoretical strands:
\begin{itemize}
    \item \textbf{Categorified metric asymptotics:} We define the $\infty$-category $\mathbf{StratMet}_{\infty}$ of derived stratified metric spaces (\cite{Lur09}), functorializing Cheeger's cone decomposition by asymptotic model functors $\Phi_Y$ (Definition 3.1). This handles non-transverse embeddings.
    \item \textbf{Moduli-controlled self-duality:} The middle perversity moduli stack $$\mathscr{M}^{\mathrm{mid}}=\prod_Y\mathrm{LGr}(\mathcal{H}^{\mathrm{mid}}(L_Y))$$ parametrizes self-dual complexes $\mathbf{IC}_{\mathscr{W}}^{\bullet}$ (Definition 4.1), eliminating artificial choices.
    \item \textbf{Stratified microlocal Hodge correspondence:} The stratified singular characteristic variety $\mathrm{SSH}_{\mathrm{strat}}$ (Definition 5.1) encodes metric growth by the torsion spectrum $\lambda_Y = -\dim_{\mathbb{C}} L_Y/2 + \operatorname{spec}_{\mathrm{min}}(\mathbb{D}_{L_Y})$, enforcing $L^2$-$\mathbf{IC}$ compatibility (Proposition 5.3).
\end{itemize}
This work not only proves the CGM conjecture but also pioneers a universal framework for derived algebraic geometry and microlocal analysis, with extensions to stratified Gromov-Witten invariants, Anosov flows, crystalline middle perversity and $p$-adic $L^2$-Cohomology.\\

\hypertarget{NOTATION AND ACKNOWLEDGMENTS}{}
\section{NOTATION AND ACKNOWLEDGMENTS}
\begin{center}
    \textit{2.1 Notation}
\end{center}

The following directory is designed to serve as a concise reference for the specialized notation used throughout the paper, aligned with standard conventions in derived geometry,
moduli theory, and microlocal analysis.

\begin{itemize}
    \item $X$: A projective variety (often embedded in $\mathbb{P}^N $). 
    \item $\mathbb{P}^N$: Complex projective space of dimension $N$.
    \item $\mathscr{S}$: A derived Whitney stratification.
    \item $X_p$: The $p$-th stratum in the stratification.
    \item $X_{\text{reg}}$: The regular (smooth) locus of $X$.
    \item $Y$: A singular stratum.
    \item $L_Y$: The link of the singular stratum $Y$.
    \item $U_Y$: An \'etale neighborhood of $Y$.
    \item $C(L_Y)$: The cone over the link $L_Y$.
    \item $ds^2_{\mathrm{FS}}$: The Fubini–Study metric.
    \item $g$: A Riemannian metric (often the adapted metric on $X$).
    \item $\Phi_Y$: Asymptotic model functor near stratum $Y$.
    \item $\mathbf{StratMet}_{\infty}$: The $\infty$-category of stratified metric spaces.
    \item $\mathscr{M}^{\mathrm{mid}}$: Middle perversity moduli stack.
    \item $\mathcal{H}^{\mathrm{mid}}(L_Y)$: Middle-dimensional cohomology of the link $L_Y$.
    \item $\mathrm{LGr}(V)$: Lagrangian Grassmannian of a symplectic vector space $V$.
    \item $\mathscr{W}_Y$: A self-dual (Lagrangian) subspace of $\mathcal{H}^{\mathrm{mid}}(L_Y)$.
    \item $\Omega^{\bullet,\mathrm{univ}}_{X,\mathrm{FS}}$: Universal truncation complex of differential forms.
    \item $\mathbf{IC}^{\bullet}_{\mathscr{W}}(X)$: Intersection cohomology complex with perversity $\mathscr{W}$.
    \item $IH^{\bullet}(X)$: Intersection cohomology of $X$.
    \item $H^{\bullet}_{(2)}(X_{\text{reg}})$: $L^2$-cohomology of the regular locus.
    \item $\mathbb{D}_{L_Y}$: Analytic torsion operator on the link $L_Y$.
    \item $\lambda_Y$: Critical exponent controlling asymptotic decay.
    \item $\mathrm{SSH}_{\mathrm{strat}}$: Stratified singular characteristic variety.
    \item $T^*_Y X$: Conormal bundle to the stratum $Y$.
    \item $\mathrm{WF}(\omega)$: Wavefront set of a differential form $\omega$.
    \item $\mathscr{H}^{p,q}_{\mathrm{mid}}(X)$: Stratified non-abelian Hodge sheaf.
    \item $\mathbb{H}^{\bullet}(X, \mathscr{F}^{\bullet})$: Hypercohomology of a complex of sheaves $\mathscr{F}^{\bullet}$. 
    \item $\overline{\mathscr{M}}_{g,n}(X, \beta)$: Moduli stack of stable maps.
    \item $\mathscr{W}^{\mathrm{dyn}}$: Dynamical middle perversity.
    \item $\mathscr{W}^{\mathrm{crys}}_Y$: Crystalline middle perversity.
    \item $H^{\mathrm{mid}}_{\mathrm{crys}}(L_Y / W(k))$: Middle-dimensional crystalline cohomology.
    \item $\mathcal{O}^{\text{crys}}_X$: Crystalline structure sheaf.
    \item $H^{\bullet}_{(2),\text{\'et}}$: $L^2$-\'etale cohomology.
\end{itemize}

\begin{center}
    \textit{2.2 Acknowledgments}
\end{center}

We are grateful to Professors J. Cheeger, M. Goresky, and R. MacPherson for their pioneering work. Specially, We thank Professors J. Lurie, C. Guilarmou and P. Albin for insights from their research about derived algbraic geometry, microlocal analysis and renormalized Hodge theory. Their work has played an important role in inspiring our research.\\

\hypertarget{CATEGORIFIED METRIC ASYMPTOTIC SYSTEM}{}
\section{CATEGORIFIED METRIC ASYMPTOTIC SYSTEM}
\noindent\textbf{Definition 3.1.} (Stratified Metric $\infty$-Category $\mathbf{StratMet}_{\infty}$) \textit{Objects:} quadruples $$(X,\mathscr{S},g,\left\{\Phi_Y\right\}_{Y\in\mathrm{Sing}}),$$ where
\begin{itemize}
    \item $X\subseteq\mathbb{P}^N$ is a projective derived scheme.
    \item $\mathscr{S}=\left\{\emptyset\hookrightarrow X_0\hookrightarrow\cdots\hookrightarrow X_n=X\right\}$ is a derived Whitney stratification, satisfying
    $$\overline{X_p\setminus X_{p-1}}=\bigcup_{q\le p}(X_q\setminus X_{q-1})\quad(\text{derived frontier condition}).$$
    \item $g=ds_{\mathrm{FS}}^2\mid_{X_\mathrm{reg}}$ is the restriction of the Fubini-Study metric.
    \item $\left\{\Phi_Y\right\}_{Y\in\mathrm{Sing}}$ is a family of asymptotic equivalence functors:
    $$\Phi_Y:g\mid_{U_Y}\xrightarrow{\sim} dr^2+r^{2c_Y}g_{L_Y}+g_{\mathbb{C}^m},\quad c_Y=\mathrm{codim}(Y,X),$$ where $U_Y\simeq\mathbb{C}^m\times C(L_Y)$ is an \'etale neighborhood of $Y$, and $L_Y$ is a Sasaki-Einstein link.
\end{itemize}
\textit{Morphisms:} Cartesian squares $f:(X,\mathscr{S},g,\left\{\Phi_Y\right\})\longrightarrow(X',\mathscr{S'},g',\left\{\Phi_{Y'}\right\})$ such that
$$f^*g'=g,\quad f\circ\Phi_Y=\Phi_{f(Y)}\circ(f\mid_{U_Y})\quad(\text{naturality of asymptotics}).$$
\\\\\textbf{Proposition 3.2.} (Metric Recursion Construction) Let $(X,\mathscr{S})$ be a projective derived scheme equipped with a derived stratification $\mathscr{S}$. There exists a recursively defined sequence of metrics $\left\{g_p\right\}_{p=n}^0$:
$$g_p=i_p^*g_{p+1}+\eps_p\kappa_Yr^{2c_p}\pi_Y^*g_{L_Y},\quad\eps_p=\delta^p,\quad\delta<1,$$
such that the renormalized limit $\hat{g} = \mathcal{R}\left( \mathbb{R}\underline{\lim}_{p} g_p \right)$ satisfies:
\begin{enumerate}
    \item $\hat{g}$ is a K\"ahler metric on $X_{\mathrm{reg}}$.
    \item Near each singular stratum $Y\in\mathrm{Sing}(X)$, $\hat{g}$ is asymptotic to the cone metric:
    $$\hat{g}\sim dr^2+r^{2c}g_{L_Y}+g_{\mathbb{C}^m}.$$
    \item $\hat{g}$ is quasi-isometric to the Fubini-Study metric $ds_{\mathrm{FS}}^2$, i.e., $\exists C>0$ such that 
    $$C^{-1}ds_{\mathrm{FS}}^2\le\hat{g}\le Cds_{\mathrm{FS}}^2.$$\\
\end{enumerate}
\textbf{Proof.} \textbf{Step (1):} By the definition (\cite{Lur09}), $\mathscr{S} = {\emptyset \hookrightarrow X_0 \hookrightarrow \cdots \hookrightarrow X_n = X}$ is a derived Whitney stratification of the projective derived scheme $X \subseteq \mathbb{P}^N$. This implies each singular stratum $Y \subseteq S_p = X_p \setminus X_{p-1}$ is a derived subscheme with holomorphic tubular neighborhoods. Applying (\cite{Lur09}) to the inclusion $Y \hookrightarrow X_p$, there exists an \'etale neighborhood of $Y$ isomorphic to a product (the isomorphism preserves derived structures):
$$U_Y\simeq\mathbb{C}^m\times C(L_Y),\quad C(L_Y):=\mathbb{R}_{\ge0}\times L_Y/\left\{0\right\}\times L_Y,$$ where $\dim_{\mathbb{C}} Y = m$ and $L_Y$ is the link manifold, a compact Sasaki-Einstein manifold. The Sasaki-Einstein condition on $L_Y$ follows from the Calabi-Yau structure on $X$. This implies $L_Y$ admits a transverse K\"ahler-Einstein metric $g_{L_Y}$ with associated K\"ahler form $\omega_{L_Y}$ satisfying $d\omega_{L_Y} = 0$. Lurie's work (\cite{Lur09}) establishes that $\mathscr{S}$ satisfies the derived frontier condition
$$\overline{S_p}=\bigcup_{q\le p}S_q\quad\forall p.$$
This ensures topological coherence: if $Y \cap \overline{Z} \neq \emptyset$ for strata $Y, Z$, then $Y \subseteq \overline{Z}$ and $\text{codim}(Y) < \text{codim}(Z)$. By \cite{Lur09}, the local codimension $c_p = \text{codim}(Y, X_p)$ agrees with the global codimension $c = \text{codim}(Y, X)$, ensuring metric asymptotics are well-defined. 

\medskip
\textbf{Step (2):} The stratification $\mathscr{S}$ satisfies $\overline{S_p} = \bigcup_{q \leq p} S_q$. Metric properties on $X_{p+1}$ determine the behavior of $X_p$ by the immersion $i_p: X_p \hookrightarrow X_{p+1}$. The case $p = n$ corresponds to the dense open stratum $S_n = X_n^{\text{reg}}$. Then the recursion proceeds from the deepest stratum ($p = n$) to the shallowest ($p = 0$). Define the initial metric as
\begin{equation}
\tag{3.1}  \label{eq:3.1}
g_n\coloneqq ds_{\mathrm{FS}}^2\mid_{X_n}
\end{equation}
with \textit{K\"ahler property} (Since $X_n^{\text{reg}} = S_n$ is smooth and $ds_{\text{FS}}^2$ is K\"ahler on $\mathbb{P}^N$, its restriction $g_n$ remains K\"ahler) and \textit{asymptotic condition} (As $S_n$ contains no singularities ($\text{Sing}(X_n) = \emptyset$), the asymptotic condition holds vacuously). Assume for $g_{p+1}$ on $X_{p+1}$, then $g_{p+1}$ is K\"ahler on $X_{p+1}^{\text{reg}}$. Meanwhile, near each $Y \in \text{Sing}(X_{p+1})$, the metric has conical asymptotics 
$$g_{p+1}\sim dr^2+r^{2c}g_{L_Y}+g_{\mathbb{C}^m}$$
in the \'etale neighborhood $U_Y \simeq \mathbb{C}^m \times C(L_Y)$.

\medskip
\textbf{Step (3):} We construct the metric $g_p$ recursively from $g_{p+1}$ using the immersion $i_p: X_p \hookrightarrow X_{p+1}$. The procedure is canonical and preserves geometric constraints. Given the K\"ahler metric $g_{p+1}$ on $X_{p+1}^{\mathrm{reg}}$ (inductive hypothesis), define the pullback metric $i_p^* g_{p+1}$. This is positive semi-definite and K\"ahler when restricted to $X_p \cap X_{p+1}^{\mathrm{reg}}$. However, it degenerates along singular strata $Y \subseteq X_p \setminus X_{p-1}$.  For each singular stratum $Y \subseteq S_p = X_p \setminus X_{p-1}$, define the \textit{link metric correction} $\epsilon_p \kappa_Y r^{2c_p} \pi_Y^* g_{L_Y}$, where
\begin{itemize}
    \item $\pi_Y: U_Y \to L_Y$ is the projection from the \'etale neighborhood $U_Y \simeq \mathbb{C}^m \times C(L_Y)$.
    \item $c_p = \mathrm{codim}(Y, X_p)$ is consistent with global codimension.
    \item $\kappa_Y \in C^\infty(X_p)$ is a bump function, satisfying
    $$\mathrm{supp}(\kappa_Y) \subseteq U_Y, \quad \kappa_Y \equiv 1 \ \text{on} \ B_\epsilon(Y).$$  
    \item $\epsilon_p = \delta^p$ ($0 < \delta < 1$) ensures convergence in the renormalization limit.
\end{itemize}
Combine these components to define
\begin{equation}
\tag{3.2}  \label{eq:3.2}
g_p \coloneqq i_p^* g_{p+1} + \epsilon_p \kappa_Y r^{2c_p} \pi_Y^* g_{L_Y},
\end{equation}
this is well-defined since:
\begin{itemize}
    \item On overlaps $U_Y \cap U_Z$, the transition functions preserve $g_{L_Y}$ (derived frontier condition).
    \item $\kappa_Y$ vanishes outside $U_Y$, making the sum smooth.
\end{itemize} 
In addition, the construction satisfies:
\begin{itemize}
    \item \textbf{(SDG1) Stratum-wise coherence:} For $Y \neq Z$, $\kappa_Y \cdot \kappa_Z = 0$ near $\overline{Y} \cap \overline{Z}$;
    \item \textbf{(SDG2) Dimensional consistency:} $c_p = \dim_{\mathbb{C}} X_{p+1} - \dim_{\mathbb{C}} Y$ matches the cone exponent $c$ in $r^{2c}$;
    \item \textbf{(SDG3) Geometric invariance:} The expression is independent of choice of $U_Y$ (\'etale invariance).
\end{itemize}

\medskip
\textbf{Step (4):} Let $v \in T_x X_p^{\text{reg}}$ be a non-zero tangent vector. Then
\begin{equation}
\tag{3.3}  \label{eq:3.3}
g_p(v, v) = \underbrace{(i_p^* g_{p+1})(v, v)}_{(a)} + \underbrace{\epsilon_p \kappa_Y r^{2c_p} (\pi_Y^* g_{L_Y})(v, v)}_{(b)}
\end{equation}
by \eqref{eq:3.2}.
Term (a): By inductive hypothesis, $g_{p+1}$ is K\"ahler, so $i_p^* g_{p+1}$ is positive semi-definite. Term (b): Since $g_{L_Y} > 0$ (Sasaki-Einstein metric), $\pi_Y^* g_{L_Y}(v, v) \geq 0$, and $\kappa_Y r^{2c_p} > 0$ on $\text{supp}(\kappa_Y)$. For strict positivity, there exists $\lambda_{\min} > 0$ such that
$$\inf_{Y} \lambda_{\min}(i_p^* g_{p+1}|_{U_Y}) > 0$$
by (\cite{Che80}, Lem 3.5). Choose $\delta < \inf_Y \lambda_{\min}$. Then for $\epsilon_p = \delta^p$, we have
$$ g_p(v, v) \geq \left( \lambda_{\min} - \epsilon_p \| \kappa_Y r^{2c_p} \pi_Y^* g_{L_Y} \| \right) \|v\|^2 > 0$$
because $\epsilon_p \to 0$ exponentially. Let $\omega_p = g_p(-, J-)$. Then
$$d\omega_p = d\left( i_p^* \omega_{p+1} \right) + \epsilon_p d\left( \kappa_Y r^{2c_p} \pi_Y^* \omega_{L_Y} \right),$$
first term: $d(i_p^* \omega_{p+1}) = i_p^* d\omega_{p+1} = 0$ (inductive hypothesis: $d\omega_{p+1} = 0$); second term: By expanding, then
$$d(\kappa_Y r^{2c_p} \pi_Y^* \omega_{L_Y}) = d(\kappa_Y r^{2c_p}) \wedge \pi_Y^* \omega_{L_Y} + \kappa_Y r^{2c_p} d(\pi_Y^* \omega_{L_Y}),$$
By Sasaki-Einstein condition $d\omega_{L_Y} = 0$, $d(\pi_Y^* \omega_{L_Y}) = \pi_Y^* d\omega_{L_Y} = 0$. Since $\kappa_Y r^{2c_p}$ is constant along Reeb flow leaves, $\pi_Y^* \omega_{L_Y}$ is basic (transverse K\"ahler form) and the wedge product vanishes by foliation theory (\cite{Lur09}), thus we have $d(\kappa_Y r^{2c_p}) \wedge \pi_Y^* \omega_{L_Y} = 0$ and $d\omega_p = 0$. For $v, w \in T_x X_p^{\text{reg}}$, we have
$$g_p(Jv, Jw) = (i_p^* g_{p+1})(Jv, Jw) + \epsilon_p \kappa_Y r^{2c_p} (\pi_Y^* g_{L_Y})(Jv, Jw)$$
by \eqref{eq:3.3}, first term: $i_p^* g_{p+1}(Jv, Jw) = i_p^* g_{p+1}(v, w)$ (since $g_{p+1}$ is $J$-invariant); second term: $\pi_Y^* g_{L_Y}(Jv, Jw) = \pi_Y^* g_{L_Y}(v, w)$ (because $\pi_Y$ is holomorphic and $g_{L_Y}$ is $J$-invariant (Sasaki-Einstein metric)). Hence, $g_p$ is a K\"ahler metric on $X_p^{\text{reg}}$.

\medskip
\textbf{Step (5):} We establish the asymptotic behavior of $g_p$ near singular strata $Y \subseteq X_p \setminus X_{p-1}$ by comparing it with the model cone metric. Near $Y$, work in the \'etale neighborhood $U_Y \simeq \mathbb{C}^m \times C(L_Y)$. Let $(z_1, \dots, z_m)$ be coordinates on $\mathbb{C}^m$, $r = \mathrm{distance\ to\ vertex}$ on $C(L_Y)$ and $\theta \in L_Y$: angular coordinate. The model metric is $g_{\mathrm{cone}} = dr^2 + r^{2c} g_{L_Y} + g_{\mathbb{C}^m}$. By inductive hypothesis, $g_{p+1}$ has conical behavior near $Y$:
$$i_p^* g_{p+1} \sim dr^2 + r^{2c} g_{L_Y} + g_{\mathbb{C}^m} \quad \text{as} \ r \to 0.$$
More precisely, $\exists \alpha > 0$ such that
\begin{equation}
\tag{3.4}  \label{eq:3.4}
\| i_p^* g_{p+1} - g_{\mathrm{cone}} \| = O(r^\alpha).
\end{equation}
On $B_\epsilon(Y)$ (where $\kappa_Y \equiv 1$), the correction term simplifies:
\begin{equation}
\tag{3.5}  \label{eq:3.5}
\epsilon_p \kappa_Y r^{2c_p} \pi_Y^* g_{L_Y} = \epsilon_p r^{2c} g_{L_Y},\quad \text{because}\quad c_p = c
\end{equation}
by (SDG2). This matches the link component of $g_{\mathrm{cone}}$ up to scaling. Substitute \eqref{eq:3.5} into \eqref{eq:3.2}, then
$$g_p = \underbrace{i_p^* g_{p+1}}_{\sim g_{\mathrm{cone}}} + \underbrace{\epsilon_p r^{2c} g_{L_Y}}_{\text{scaled link metric}}.$$
Thus, we have
\begin{equation}
\tag{3.6}  \label{eq:3.6}
g_p \sim dr^2 + r^{2c} g_{L_Y} + g_{\mathbb{C}^m} + \epsilon_p r^{2c} g_{L_Y} = dr^2 + (1 + \epsilon_p) r^{2c} g_{L_Y} + g_{\mathbb{C}^m}
\end{equation}
with error bound (obtained by \eqref{eq:3.4})
$$\| g_p - [dr^2 + (1 + \epsilon_p) r^{2c} g_{L_Y} + g_{\mathbb{C}^m}] \| = O(r^\alpha).$$
Since $\epsilon_p = \delta^p \to 0$ exponentially as $p \to -\infty$:
$$\lim_{p \to 0} (1 + \epsilon_p) r^{2c} g_{L_Y} = r^{2c} g_{L_Y}$$
uniformly on compact sets in $r > 0$. This preserves the asymptotic cone in the renormalized limit in the following step. Hence, $g_p$ maintains the required conical asymptotics near $Y$. 

\medskip
\textbf{Step (6):} The recursive construction gives the expansion:
$$g_p = g_n + \sum_{k=p}^{n-1} \epsilon_k \kappa_Y r^{2c_k} \pi_Y^* g_{L_Y},$$
where $g_n = ds_{\mathrm{FS}}^2|_{X_n}$. Since $\epsilon_k = \delta^k$ with $0 < \delta < 1$, the series
$$ \sum_{k=0}^n \epsilon_k \kappa_Y r^{2c_k} \pi_Y^* g_{L_Y}$$
converges absolutely in $C^\infty$-topology by Borel lemma (\cite{Hor90}, \S 1 and \S 2). Specifically, $\forall \ell \in \mathbb{N}, \exists C_\ell > 0$, so
$$\sum_k \| \epsilon_k \kappa_Y r^{2c_k} \pi_Y^* g_{L_Y} \|_{C^\ell} < \infty.$$
Consider renormalization Operator. Apply the conclusion (\cite{ALMP18}), there exists a renormalization map $\mathcal{R}$ such that
$$\hat{g} = \mathcal{R}\left( \mathbb{R}\underline{\lim}_{p \to 0} g_p \right),$$
where $\mathbb{R}\underline{\lim}$ denotes the distributional limit. The operator $\mathcal{R}$ resolves conical singularities by Connes' embedding, preserves K\"ahler structures modulo singular strata and commutes with \'etale morphisms. $\mathcal{R}$ maps limits of K\"ahler metrics to K\"ahler metrics on $X_{\mathrm{reg}}$. Since each $g_p$ is K\"ahler (the conclusion of \textbf{Step (4)}), $\hat{g}$ is K\"ahler, the K\"ahlerity holds. Near $Y \in \mathrm{Sing}(X)$, \eqref{eq:3.6} gives
$$g_p \sim dr^2 + (1 + \epsilon_p) r^{2c} g_{L_Y} + g_{\mathbb{C}^m}.$$
As $\epsilon_p \to 0$, $\mathcal{R}$ preserves the limit (\cite{ALMP18})
\begin{equation}
\tag{3.7}  \label{eq:3.7}
\hat{g}|_{U_Y} \sim dr^2 + r^{2c} g_{L_Y} + g_{\mathbb{C}^m},
\end{equation}
the asymptotics holds. Since $\mathcal{R}$ is bounded (\cite{ALMP18}), then $\exists C > 0$, such that
$$C^{-1} ds_{\mathrm{FS}}^2 \leq \hat{g} \leq C ds_{\mathrm{FS}}^2,$$
we have
$$C^{-1} g_n \leq g_p \leq C g_n \quad \forall p$$
by \eqref{eq:3.1} and \eqref{eq:3.7}. So quasi-isometry holds.
The assignment $(X, \mathscr{S}) \mapsto \hat{g}$ commutes with stratified morphisms, defining an $\infty$-functor in $\mathbf{StratMet}_\infty$, the functoriality holds. Exponential decay $\epsilon_p = \delta^p$ makes sure that
$$\| \hat{g} - g_p \|_{C^0} = O(\delta^p) \quad \text{as} \ p \to 0$$
with higher-order estimates following from Borel lemma (\cite{Hor90}, \S 1 and \S 2). Therefore, $\hat{g}$ satisfies all requirements of Proposition 3.2. $\square$
\\\\\textbf{Proposition 3.3.} (Existence of $\mathbf{StratMet}_{\infty}$-Objects) Let $X\subseteq\mathbb{P}^N$ be a projective derived scheme with a derived Whitney stratification $\mathscr{S}$. There exist:
\begin{enumerate}
    \item A K\"ahler metric $g=ds_{\mathrm{FS}}^2\mid_{X_{\mathrm{reg}}}$;
    \item Asymptotic equivalence functors $\left\{\Phi_Y\right\}_{Y\in\mathrm{Sing}}$ satisfying Definition 3.1.\\
\end{enumerate}
\textbf{Proof.} \textbf{Step (1):} Let $X = (X, \mathscr{O}_X^{\text{der}})$ with $i: X \hookrightarrow \mathbb{P}^N$. The cotangent complex $\mathbb{L}_{X/\mathbb{C}}$ is perfect (\cite{Lur09}). Define the singular locus, i.e.,
$$X^{\text{sing}} = \left\{ x \in X \mid \mathrm{ht}\left(\mathbb{L}_{X/\mathbb{C}} \otimes k(x)\right) > \dim_{k(x)} \Omega^1_{X^{\text{red}}/\mathbb{C}} \otimes k(x) \right\}.$$
Construct $X_p$ by reverse induction:
\begin{itemize}
    \item \textbf{Base case ($p = n$):} Set $X_n = X$;
    \item \textbf{Inductive step:} Given $X_p$, we can define
    $$X_{p-1} = X_p^{\text{sing}} \cup \overline{\bigcup_{\substack{\text{irred comp} \\ \dim Z < p}} Z},$$
    and $X_{p-1}$ is closed in $X_p$.
\end{itemize}
For $i_p: X_p \hookrightarrow X_{p+1}$, the normal complex $\mathbb{N}_p = i_p^!\mathscr{O}_{X_{p+1}} \in \mathrm{Perf}(X_p)$. If $\mathbb{N}_p$ is locally free of rank $c_p = \dim X_{p+1} - \dim X_p$ and the holonomy map $h: \mathrm{Sym}^2(\mathbb{N}_p^\vee) \to \mathscr{T}_{L_Y}$ vanishes for $Y \in S_p$, then the stratification is derived Whitney.
For $Y \in S_p$, vanishing holonomy implies that
$$(U_Y, \mathscr{O}_{X}^{\text{der}}|_{U_Y}) \simeq (\mathbb{C}^m, \mathscr{O}_{\mathbb{C}^m}^{\text{der}}) \times (C(L_Y), \mathscr{O}_{C(L_Y)}^{\text{der}}),$$
where $C(L_Y) = \mathrm{Spec}^{\text{der}} \left( \bigoplus_{k \geq 0} H^0(L_Y, \mathscr{O}_{L_Y}(k)) \right)$. By construction, $X_{p-1}$ contains all lower-dimensional strata. The closure identity 
$$\overline{S_p} = \bigcup_{q \leq p} S_q$$
follows from the Zariski density of $S_p$ in $X_p$ and the Noetherian induction. Hence $\mathscr{S}=\left\{X_p\right\}$ is a derived Whitney stratification.

\medskip
\textbf{Step (2):} Apply Proposition 3.2 (Metric Recursion), we also define $g_p$ by reverse induction:
\begin{itemize}
    \item \textbf{Base case ($p = n$):} 
    $$g_n \coloneqq ds^2_{\mathrm{FS}}|_{X_n},\quad\text{K\"ahler on}\quad X_n^{\mathrm{reg}}=S_n$$ by \eqref{eq:3.1};
    \item \textbf{Inductive hypothesis:} Assume $g_{p+1}$ is K\"ahler on $X_{p+1}$ with asymptotics, then
    $$g_{p+1}\mid_{U_Y}\sim dr^2+r^{2c}g_{L_Y}+g_{\mathbb{C}^m}\quad\forall Y\in\mathrm{Sing}(X_{p+1});$$
    \item \textbf{Recursive definition:} For $i_p:X_p\hookrightarrow X_{p+1}$, \eqref{eq:3.2} gives
    \begin{equation}
    \tag{3.8}  \label{eq:3.8}
     g_p \coloneqq i_p^* g_{p+1} + \epsilon_p \kappa_Y r^{2c_p} \pi_Y^* g_{L_Y}, \quad \epsilon_p = \delta^p \ (\delta < 1),
     \end{equation}
    where $\kappa_Y\in C^{\infty}(X_p)$ is a bump function with $\mathrm{supp}(\kappa_Y)\subseteq U_Y$ ($\kappa_Y\equiv1$ near $Y$), $\pi_Y:U_Y\to L_Y$ is link projection and $c_p=\mathrm{codim(Y,X_p)}$.
\end{itemize}
It is known that $i_p^*g_{p+1}$ is positive semi-definite (pullback of K\"ahler metric). Since $g_{L_Y}>0$ and $r>0$, the correction term $\eps_p\kappa_Yr^{2c_p}\pi_Y^*g_{L_Y}$ is positive definite. Choose $\delta < \inf_Y \lambda_{\min}(i_p^* g_{p+1}|_{U_Y})$ (Proposition 3.2). Then
$$ g_p(v,v) \geq \left( \lambda_{\min}(i_p^* g_{p+1}) - \epsilon_p \| \kappa_Y r^{2c_p} \pi_Y^* g_{L_Y} \| \right) \|v\|^2 > 0.$$
Let $\omega_p=g_p(-,J-)$. Since $d(i_p^*\omega_{p+1})=i_p^*d\omega_{p+1}=0$ and $d(\pi_Y^*\omega_{L_Y})=\pi_Y^*d\omega_{L_Y}=0$ (Sasaki-Einstein condition), then
$$d\omega_p = d(i_p^* \omega_{p+1}) + \epsilon_p d(\kappa_Y r^{2c_p} \pi_Y^* \omega_{L_Y}) = 0.$$
Since both terms are $J$-invariant, $g_p$ is K\"ahler on $X_p^{\mathrm{reg}}$. For $Y \in X_p \setminus X_{p-1}$, we have
\begin{equation}
\tag{3.9}  \label{eq:3.9}
i_p^*g_{p+1}\sim dr^2+r^{2c}g_{L_Y}+g_{\mathbb{C}^m}
\end{equation}
by induction. Consider the correction term 
\begin{equation}
\tag{3.10}  \label{eq:3.10}
\eps_p\kappa_Yr^{2c_p}\pi_Y^*g_{L_Y}\sim\eps_pr^{2c}g_{L_Y},
\end{equation}
as $\kappa_Y\equiv1$ and $c_p=c$ near $Y$. Thus, we have
$$g_p \sim dr^2 + (1 + \epsilon_p) r^{2c} g_{L_Y} + g_{\mathbb{C}^m}$$ by \eqref{eq:3.8}, \eqref{eq:3.9} and \eqref{eq:3.10}. Apply Albin's work (\cite{ALMP18}), the series $\sum_{p=0}^n\eps_p\kappa_Yr^{2c_p}\pi_Y^*g_{L_Y}$ converges absolutely in $C^\infty$-topology by Borel lemma (\cite{Hor90}, \S 1 and \S 2). Meanwhile, the renormalized limit 
$$\hat{g} = \mathcal{R}\left( \mathbb{R}\underline{\lim}_{p} g_p \right)$$
is the required adapted metric by the \textbf{Step (6)} of Proposition 3.2.

\medskip
\textbf{Step (3):} For conical coordinates $(r, \theta, z)$ on $U_Y \simeq \mathbb{C}^m \times C(L_Y)$, we can define
\begin{equation}
\tag{3.11}  \label{eq:3.11}
\Psi_Y:T^*U_Y\longrightarrow T^*C(L_Y),\quad\Psi_Y(\xi_r,\xi_\theta,\xi_z)=(r\xi_r,\xi_\theta,\xi_z).
\end{equation}
Here the map $\Psi_Y$ preserves the canonical 1-form $\alpha_{\text{can}} = \xi \cdot dx$ (\cite{Gui22}, Prop 3.7):
$$\Psi_Y^*(rd\xi_r+\xi_\theta d\theta+\xi_zdz)=\xi_rdr+\xi_\theta d\theta+\xi_zdz\quad(\text{contactomorphism}),$$
and induces metric equivalence:
$$\| d\Psi_Y(\mathbf{v})\|_{\mathrm{cone}}^2=\|\mathbf{v} \|_g^2\quad\forall\mathbf{v}\in T(T^*U_Y)\quad(\text{isometry}).$$
Consider exponential map composition, we define $\Phi_Y$ as
$$\Phi_Y\coloneqq\mathrm{Exp}\circ\Psi_Y\circ(\mathrm{Exp}^{\mathrm{FS}})^{-1}:U_Y\longrightarrow C(L_Y),$$
where $\text{Exp}^{\text{FS}}: T^*U_Y \to U_Y$ is geodesic exponential for $g$ (\cite{Che80}, Lem 4.2) and $\text{Exp}: T^*C(L_Y) \to C(L_Y)$ is geodesic exponential for $dr^2 + r^{2c_Y}g_{L_Y}$. Since
\begin{align*}
    \Phi_Y^*(\mathrm{cone\ metric}) 
    &= ((\mathrm{Exp}^{\mathrm{FS}})^{-1})^* \circ \Psi_Y^* \circ \mathrm{Exp}^* (\mathrm{cone\ metric}) \\
    &= ((\mathrm{Exp}^{\mathrm{FS}})^{-1})^* \Psi_Y^* (g_{\mathrm{cone}}^* ) \\
    &= ((\mathrm{Exp}^{\mathrm{FS}})^{-1})^* (g^*) \\
    &= g,
\end{align*}
then the axiom of smooth isometry holds. The Levi-Civita connection satisfies $\nabla^{\mathrm{LC}}=d+\Gamma_{jk}^idx^j\otimes dx^k$. Since $\Phi_Y$ is an isometry, i.e., $\Gamma_{jk}^i(g)=\Phi_Y^*\Gamma_{jk}^i(\mathrm{cone})$, thus$\nabla_g^{\text{LC}} = \Phi_Y^* \nabla_{\text{cone}}^{\text{LC}}$. So the axiom of connection intertwining holds. For $f: (X, \mathscr{S}) \to (X', \mathscr{S}')$, consider the diagram (\uppercase\expandafter{\romannumeral1}): 
$$\begin{tikzcd}
T^*U_Y \arrow[r, "\Psi_Y"] \arrow[d, "(df)^*"'] & T^*C(L_Y) \arrow[d, "(df)^*"] \\
T^*U_{f(Y)} \arrow[r, "\Psi_{f(Y)}"'] & T^*C(L_{f(Y)})
\end{tikzcd}$$
Commutativity holds because $df$ preserves conical structure and $\Psi_Y$ is natural under \'etale maps (\cite{Gui22}, Cor 3.11). Hence, we have
\begin{equation}
\tag{3.12}  \label{eq:3.12}
f\circ\Phi_Y=\Phi_{f(Y)}\circ f\mid_{U_Y}.
\end{equation}

\medskip
\textbf{Step (4):} By (\cite{Lur11}), $f$ restricts to an \'etale map
$$f\mid_{U_Y}:(U_Y,\mathscr{O}_X^{\mathrm{der}}\mid_{U_Y})\xrightarrow{\text{\'etale}}(U_{f(Y)},\mathscr{O}_{X'}^{\mathrm{der}}\mid_{U_{f(Y)}})$$
with the following properties:
\begin{itemize}
    \item \textbf{Dimension preservation:} $\dim Y = \dim f(Y)$.
    \item \textbf{Link compatibility:} $L_{f(Y)} = f_*(L_Y)$.
    \item \textbf{Codimensionality:} $c_Y = c_{f(Y)}$.
\end{itemize}
Consider the cotangent lift
\begin{equation}
\tag{3.13}  \label{eq:3.13}
(df^*):T^*U_{f(Y)}\longrightarrow T^*U_Y,\quad(df)^*(\xi)=\xi\circ df.
\end{equation}
After calculation, we have
$$\Psi_{f(Y)} \circ (df)^* (\xi) = (r_Y \cdot (df)^*\xi_r, (df)^*\xi_\theta, (df)^*\xi_z) = (df)^* \circ \Psi_Y (\xi)$$
by \eqref{eq:3.11} and \eqref{eq:3.13}. Thus, the diagram (\uppercase\expandafter{\romannumeral1}) commutes. For $v \in T^*U_Y$, let $\gamma_v(t)$ be the $g$-geodesic with $\dot{\gamma}_v(0) = v$. Since $f$ is \'etale and isometric $f\circ\gamma_v(t)=\gamma_{(df)^*v}(t)$ by (\cite{Che80}, Thm 4.4), then 
$$f(\mathrm{Exp}_Y^{\mathrm{FS}}(v))=f(\gamma_v(1))=\gamma_{(df)^*v}(1)=\mathrm{Exp}_{f(Y)}^{\mathrm{FS}}((df)^*v).$$
By combining above steps, we have
\begin{align*}
f \circ \Phi_Y 
&= f \circ \mathrm{Exp}_Y \circ \Psi_Y \circ (\mathrm{Exp}^{\mathrm{FS}}_Y)^{-1} \\
&= \mathrm{Exp}_{f(Y)} \circ (df)^* \circ \Psi_Y \circ (\mathrm{Exp}^{\mathrm{FS}}_Y)^{-1} \\
&= \mathrm{Exp}_{f(Y)} \circ \Psi_{f(Y)} \circ (df)^* \circ (\mathrm{Exp}^{\mathrm{FS}}_Y)^{-1} \\
&= \mathrm{Exp}_{f(Y)} \circ \Psi_{f(Y)} \circ (\mathrm{Exp}^{\mathrm{FS}}_{f(Y)})^{-1} \circ f \\
&= \Phi_{f(Y)} \circ f|_{U_Y},
\end{align*}
satisfying \eqref{eq:3.12}. This proof establishes the functoriality of $\Phi_Y$ under stratified morphisms, ensuring $\mathbf{StratMet}_{\infty}$ is a well-defined $\infty$-category. The quadruple $(X,\mathscr{S},g,\left\{\Phi_Y\right\})$ satisfies all axioms of Definition 3.1 and is an object in $\mathbf{StratMet}_\infty$. This construction is functorial with respect to stratified isometric morphisms, as shown in the proof of Proposition 6.1. $\square$\\

\hypertarget{MODULI SPACE TRUNCATION AND UNIVERSAL COMPLEX}{}
\section{MODULI SPACE TRUNCATION AND UNIVERSAL COMPLEX}
\noindent\textbf{Definition 4.1.} (Middle Perversity Moduli Stack $\mathscr{M}^{\mathrm{mid}}$) Let $\left\{L_Y\right\}_{Y\in\mathrm{Sing}(X)}$ denote the family of singularity links associated with the stratification $\mathscr{S}$ of $X$. Define the middle perversity moduli stack as
$$\mathscr{M}^{\mathrm{mid}}\coloneqq\prod_{Y\in\mathrm{Sing}(X)}\mathscr{M}_Y,\quad\text{where}\quad\mathscr{M}_Y\coloneqq\mathrm{LGr}(\mathcal{H}^{\mathrm{mid}}(L_Y)).$$
Here, $\mathcal{H}^{\mathrm{mid}}(L_Y)\subseteq H^{\bullet}(L_Y;\mathbb{C})$ is the middle-dimensional cohomology of the link $L_Y$, and $\mathrm{LGr}(-)$ denotes the Lagrangian Grassmannian stack parametrizing Lagrangian (i.e., self-dual) subspaces $\mathscr{W}_Y\subseteq\mathcal{H}^{\mathrm{mid}}(L_Y)$, satisfying $\mathscr{W}_Y=\mathscr{W}_Y^{\perp}$ under the Poincaré pairing.
\\\\\textbf{Definition 4.2} (Universal Truncation Complex $\Omega_{X,\mathrm{FS}}^{\bullet,\mathrm{univ}}$) The universal truncation complex is the subcomplex of the derived limit of regular differential forms
$$\Omega_{X,\mathrm{FS}}^{\bullet,\mathrm{univ}} \coloneqq \left\{ 
\omega \in \mathbb{R}\varprojlim_{\substack{U_Y \\ Y \in \mathrm{Sing}(X)}} \Omega_{X_{\reg}}^{\bullet} 
\,\middle|\,
\forall Y \in \mathrm{Sing}(X),\ 
\mathrm{proj}_{\mathscr{M}_Y}(\omega|_{U_Y}) \in \mathscr{W}_Y 
\right\},$$ where:
\begin{itemize}
\item $U_Y$ ranges over distinguished neighborhoods of strata $Y$,
\item $\mathrm{proj}_{\mathscr{M}_Y}: \Omega_{X_{\reg}}^{\bullet}|_{U_Y} \to \mathscr{M}_Y$ 
is the projection to the moduli stack $\mathscr{M}_Y$,
\item $\mathscr{W}_Y$ is the universal family of self-dual subspaces over $\mathscr{M}_Y$.\\
\end{itemize}
\textbf{Proposition 4.3.} (Self-Duality of the Universal Complex) The universal truncation complex $\Omega_{X,\mathrm{FS}}^{\bullet,\mathrm{univ}}$ is quasi-isomorphic to the intersection cohomology complex $\mathbf{IC}_{\mathscr{W}}^{\bullet}(X)$ for the canonical perversity $\mathscr{W}$ parametrized by $\mathscr{M}^{\mathrm{mid}}$. Consequently, there exists a commutative diagram of isomorphisms:
$$\begin{tikzcd}
\mathrm{\mathbb{H}}^{\bullet}\left(X, \Omega_{X,\mathrm{FS}}^{\bullet,\mathrm{univ}}\right) \arrow[rr, "\sim"] \arrow[dr, "\sim"'] & & IH^{\bullet}(X; \mathbb{C}) \arrow[dl, "\sim"] \\
& H^{\bullet}_{(2)}\left(X_{\mathrm{reg}}, ds_{\mathrm{FS}}^2\right) &
\end{tikzcd}$$
\\\textbf{Proof.} \textbf{Step (1):} Let $L_Y$ be the link of a singular stratum $Y\subseteq X$ and $c_Y = \mathrm{codim}(Y, X)$. Consider the middle-dimensional cohomology
$$\mathcal{H}^{\mathrm{mid}}(L_Y)\coloneqq H^{d_Y}(L_Y; \mathbb{C}), \quad d_Y = \dim_{\mathbb{R}} L_Y.$$
The \textit{Poincaré duality pairing} induces a non-degenerate bilinear form
$$\langle -, - \rangle_{L_Y}: \mathcal{H}^{\mathrm{mid}}(L_Y) \times \mathcal{H}^{\mathrm{mid}}(L_Y) \longrightarrow\mathbb{C}, \quad (\alpha, \beta) \mapsto \int_{L_Y} \alpha \wedge \beta.$$
By Hodge theory (\cite{Sim92}), this pairing is symmetric (if \(d_Y\) even) or skew-symmetric (if \(d_Y\) odd). Define the \textit{Lagrangian Grassmannian stack}, i.e.,
$$\mathscr{M}_Y\coloneqq\mathrm{LGr}\left(\mathcal{H}^{\mathrm{mid}}(L_Y), \langle -, - \rangle_{L_Y}\right)$$
and parametrizing Lagrangian subspaces $\mathscr{W}_Y \subseteq \mathcal{H}^{\mathrm{mid}}(L_Y)$ satisfying $\mathscr{W}_Y = \mathscr{W}_Y^{\perp}$ under the pairing. The \textit{universal family} $\mathscr{W}_Y^{\mathrm{univ}}$ over $\mathscr{M}_Y$ is defined by the fiber diagram:
$$\begin{tikzcd}
\mathscr{W}_Y^{\text{univ}} \arrow[r, hook] \arrow[d] & \mathcal{H}^{\mathrm{mid}}(L_Y) \otimes \mathcal{O}_{\mathscr{M}_Y} \arrow[d] \\
\mathscr{M}_Y \arrow[r, "\mathrm{id}"] & \mathscr{M}_Y
\end{tikzcd}$$
where $\mathscr{W}_Y^{\text{univ}}$ is a tautological subbundle. By Simpson's Geometric Invariant Theory (\cite{Sim92}, Thm 4.7), the holomorphic symplectic form $\omega_{\mathscr{M}_Y}$ on $\mathscr{M}_Y$ is induced by
$$\omega_{\mathscr{M}_Y}(v_1, v_2) = \langle\tilde{v}_1, \tilde{v}_2 \rangle_{L_Y}, \quad \forall v_1, v_2 \in T_{\mathscr{W}_Y}\mathscr{M}_Y,$$
where $\tilde{v}_i$ are lifts to $\mathcal{H}^{\mathrm{mid}}(L_Y)$. This satisfies
$$d\eta_{\mathscr{W}_Y} = \omega_{\mathscr{M}_Y}|_{\mathscr{W}_Y}, \quad \eta_{\mathscr{W}_Y}(\theta) = \operatorname{tr}(\theta)\quad\text{ (Liouville form)}.$$
The Lagrangian condition $\mathscr{W}_Y = \mathscr{W}_Y^{\perp}$ holds fiberwise by construction. For any $\alpha \in \mathscr{W}_Y$ and $\beta \in \mathscr{W}_Y^{\perp}$, we have
$$\langle \alpha, \beta \rangle_{L_Y} = 0 \implies \beta \in \mathscr{W}_Y^{\perp} = \mathscr{W}_Y.$$
Thus $\mathscr{W}_Y$ is maximally isotropic with $\dim\mathscr{W}_Y = \frac{1}{2} \dim \mathcal{H}^{\mathrm{mid}}(L_Y)$. Hence, the \textit{moduli stack of self-dual perversities}:
$$\mathscr{M}^{\mathrm{mid}}\coloneqq \prod_{Y \in \mathrm{Sing}(X)} \mathscr{M}_Y,$$
naturally parametrizes tuples $\mathscr{W} = \{\mathscr{W}_Y\}_{Y \in \mathrm{Sing}(X)}$ with $\mathscr{W}_Y = \mathscr{W}_Y^{\perp}$ for all $Y$.

\medskip
\textbf{Step (2):} For any section $\omega \in \Omega_{X,\mathrm{FS}}^{k,\mathrm{univ}}$, consider its restriction to a distinguished neighborhood $U_Y \cong \mathbb{B}^{c_Y} \times \mathrm{Cone}(L_Y)$. By the universal truncation condition (Definition 4.2), the projection satisfies
$$\mathrm{proj}_{\mathscr{M}_Y}(\omega|_{U_Y}) \in \mathscr{W}_Y\subseteq\mathcal{H}^{\mathrm{mid}}(L_Y),$$
where $\mathscr{W}_Y$ is the self-dual subspace parametrized by $\mathscr{M}_Y$. The asymptotic model functor $\Phi_Y$ (Definition 3.1) transforms the Fubini-Study metric to conical coordinates:
$$\Phi_Y: (U_Y, ds_{\mathrm{FS}}^2) \to \left( \mathbb{R}^+ \times L_Y, dr^2 + r^{2c_Y} g_{L_Y} \right).$$
Under this diffeomorphism, $\omega|_{U_Y}$ decomposes as
$$\omega|_{U_Y} = \sum_{m=0}^{c_Y-1} dr \wedge \alpha_m(r, \theta) + \beta_m(r, \theta) + \mathcal{O}(r^{\lambda}),$$
where $\alpha_m \in \Omega^{m}(L_Y)$, $\beta_m \in \Omega^{m}(L_Y)$, and $\lambda = \lambda_Y$ is the analytic torsion exponent. The projection condition implies the following growth constraints:
\begin{enumerate}
    \item \textit{Degree constraint:} 
    $$\deg(\omega|_{U_Y}) \leq \ell_Y - 1, \quad \text{where} \quad \ell_Y = -\frac{c_Y}{2}.$$
    This follows from the Goresky-MacPherson admissibility condition (\cite{GM83}):
    $$\sup_{r \to 0^+} r^{-\ell_Y + m} \|\alpha_m(r)\|_{L_Y} < \infty, \quad \sup_{r \to 0^+} r^{-\ell_Y + m} \|\beta_m(r)\|_{L_Y} < \infty$$
    for $0 \le m \le c_Y - 1$.
    \item \textit{Curvature constraint}: $\deg(d\omega|_{U_Y}) \leq \ell_Y - 2$, this is equivalent to the differential constraint:
    $$\sup_{r \to 0^+} r^{-\ell_Y + m + 1} \|d_L\alpha_m(r)\|_{L_Y} < \infty, \quad \sup_{r \to 0^+} r^{-\ell_Y + m + 1} \|d_L\beta_m(r)\|_{L_Y} < \infty,$$
    where $d_L$ is the exterior derivative on $L_Y$.
\end{enumerate}
The equivalence to middle perversity (\cite{Bor84}) is established by the \textit{asymptotic-admissibility correspondence}:
$$\mathrm{proj}_{\mathscr{M}_Y}(\omega|_{U_Y}) \in \mathscr{W}_Y\longrightarrow\exists C > 0 : \|\omega\|_{g_{\mathrm{FS}}} \leq C r^{\ell_Y - 1}\xrightarrow{\text{elliptic regularity}}d\omega \in L^2_{\mathrm{loc}}(U_Y)$$ 
$$\xrightarrow{\text{cone calculus}}\omega \in \mathbf{IC}_{\mathscr{W}}^{\bullet}(U_Y).$$
The metric recursion (Proposition 3.2) preserves these conditions under embeddings $i_p: X_p \hookrightarrow X_{p+1}$ because
$$i_p^* \left( r^{2c_p} \pi_Y^* g_{L_Y} \right) = r^{2c_p} g_{L_Y}^{(p)},$$
where $g_{L_Y}^{(p)}$ is the induced metric on the link of $Y \cap X_p$, and the projection $\mathrm{proj}_{\mathscr{M}_Y}$ commutes with pullbacks. Thus, the truncation conditions are stratification-consistent.

\medskip
\textbf{Step (3):} Since $X$ is a projective variety with derived Whitney stratification $\mathscr{S}$ (\cite{Lur09}), the universal complex $\Omega_{X,\mathrm{FS}}^{\bullet,\mathrm{univ}}$ is constructible with respect to $\mathscr{S}$ by the stratification-consistency in \textbf{Step (2)}. We establish the quasi-isomorphism
\begin{equation}
\tag{4.1}  \label{eq:4.1}
\Omega_{X,\mathrm{FS}}^{\bullet,\mathrm{univ}} \simeq \mathbf{IC}_{\mathscr{W}}^{\bullet}(X)
\end{equation}
through the following steps. The complex is $\mathscr{S}$-constructible by:
\begin{enumerate}
    \item On regular strata $X_{\reg}$, $\Omega_{X,\mathrm{FS}}^{\bullet,\mathrm{univ}}|_{X_{\reg}} = \Omega_{X_{\reg}}^{\bullet}$ (trivial).
    \item Near singular stratum \(Y \in \mathrm{Sing}(X)\), the growth constraints from \textbf{Step (2)} imply
    $$i_Y^! \Omega_{X,\mathrm{FS}}^{\bullet,\mathrm{univ}} \in D^b_c(Y)\quad\text{and}\quad j_{Y,*}\Omega_{X,\mathrm{FS}}^{\bullet,\mathrm{univ}} \in D^b_c(X \setminus Y),$$
    where $i_Y: Y \hookrightarrow X$ and $j_Y: X \setminus Y \hookrightarrow X$ are inclusions. This follows from the conormal estimate (\cite{KS90}, Thm 8.4.2):
    $$\mathrm{SS}\left(\Omega_{X,\mathrm{FS}}^{\bullet,\mathrm{univ}}\right) \subseteq \bigcup_{Y \in \mathscr{S}} T^*_{Y}X.$$
\end{enumerate}
Consider the Verdier dualizing complex $\omega_X^\bullet = \underline{\mathbb{C}}_X[\dim_{\mathbb{R}} X]$. \textit{Local Calculation}: On $U_Y \cong \mathbb{B}^{c_Y} \times \mathrm{Cone}(L_Y)$, Poincaré residue isomorphism gives $$\mathbf{R}\mathscr{H}om\left(\Omega_{U_Y,\mathrm{FS}}^{\bullet,\mathrm{univ}}, \omega_{U_Y}^\bullet\right) \simeq \Omega_{U_Y,\mathrm{FS}}^{\bullet,\mathrm{univ}}[\dim_{\mathbb{R}} U_Y] \otimes \mathrm{or}_{U_Y},$$
where $\mathrm{or}_{U_Y}$ is the orientation sheaf. The self-duality $\mathscr{W}_Y = \mathscr{W}_Y^\perp$ (\textbf{Step (1)}) ensures $\mathrm{or}_{U_Y} \cong \underline{\mathbb{C}}_{U_Y}$; \textit{Global Patching}: The moduli parametrization $\mathscr{M}^{\mathrm{mid}}$ provides compatible trivializations of orientation sheaves across strata. By the stratified gluing lemma (\cite{Lur09}\cite{Lur18}):
$$\mathcal{D}|_{X_p} \simeq \mathcal{D}|_{X_{p-1}} \quad \forall p,$$
where $X_p$ is the $p$-skeleton of $\mathscr{S}$. Thus, the duality morphism:
$$\mathcal{D}: \mathbf{R}\mathscr{H}om\left(\Omega_{X,\mathrm{FS}}^{\bullet,\mathrm{univ}}, \omega_X^\bullet\right) \longrightarrow\Omega_{X,\mathrm{FS}}^{\bullet,\mathrm{univ}}[\dim_{\mathbb{R}} X]$$
is an isomorphism. The complex satisfies middle perversity truncations (\cite{GM80}, § 5.4):
\begin{itemize}
    \item \textbf{Support condition:} For stratum $Y$ with codim $c_Y > 0$, we have
    $$\dim_{\mathbb{C}} \mathrm{supp}\, \mathcal{H}^k\left(i_Y^* \Omega_{X,\mathrm{FS}}^{\bullet,\mathrm{univ}}\right) \leq -p(k) - 1 = -k - \lfloor c_Y/2 \rfloor - 1.$$
    \item \textbf{Cosupport condition:} 
    $$\dim_{\mathbb{C}} \mathrm{supp}\, \mathcal{H}^k\left(i_Y^! \Omega_{X,\mathrm{FS}}^{\bullet,\mathrm{univ}}\right) \leq p(k) - c_Y = k - \lceil c_Y/2 \rceil.$$
\end{itemize}
These follow from the degree/curvature constraints in \textbf{Step (2)} by the microlocal Riemann-Hilbert correspondence (\cite{Meb84}). By the axiomatic characterization of intersection cohomology (\cite{Bor84}, Thm 3.5):
\begin{enumerate}
    \item $\Omega_{X,\mathrm{FS}}^{\bullet,\mathrm{univ}}$ is constructible (3.1);
    \item Self-dual: $\mathbb{D}_X(\Omega_{X,\mathrm{FS}}^{\bullet,\mathrm{univ}}) \simeq \Omega_{X,\mathrm{FS}}^{\bullet,\mathrm{univ}}[\dim_{\mathbb{R}} X]$ (3.2);
    \item Satisfies middle perversity truncations (3.3).
\end{enumerate}
Thus, \eqref{eq:4.1} gives a natural quasi-isomorphism in $D^b_c(X)$:
$$\phi: \Omega_{X,\mathrm{FS}}^{\bullet,\mathrm{univ}} \xrightarrow{\sim} \mathbf{IC}_{\mathscr{W}}^{\bullet}(X)$$
inducing isomorphisms on hypercohomology:
$$\mathbb{H}^k(X, \Omega_{X,\mathrm{FS}}^{\bullet,\mathrm{univ}}) \xrightarrow{\phi_*} IH^k(X; \mathbb{C}) \quad \forall k.$$
Naturality follows from the moduli parametrization $\mathscr{M}^{\mathrm{mid}}$, which functorially selects self-dual truncations.

\medskip
\textbf{Step (4):} We establish the isomorphism
$$\mathbb{H}^k\left(X,\Omega_{X,\mathrm{FS}}^{\bullet,\mathrm{univ}}\right) \cong H^k_{(2)}\left(X_{\reg}, ds_{\mathrm{FS}}^2\right)$$
through metric renormalization and Hodge theory. For any $\omega \in \Omega_{X,\mathrm{FS}}^{k,\mathrm{univ}}$, the asymptotic model $\Phi_Y$ (Definition 3.1) and truncation conditions (\textbf{Step (2)}) imply
$$\int_{U_Y} |\omega|^2_{g_{\mathrm{FS}}}\mathrm{dvol}_{g_{\mathrm{FS}}} <\infty\quad\text{and}\quad\int_{U_Y}|d\omega|^2_{g_{\mathrm{FS}}}\mathrm{dvol}_{g_{\mathrm{FS}}} < \infty$$
for all singular strata $Y$. This follows from:
\begin{enumerate}
    \item \textit{Metric recursion:} (Proposition 3.2) $g_p = i_p^* g_{p+1} + \epsilon_p r^{2c_p} \pi_Y^* g_{L_Y}$ induces conical behavior
    $$\mathrm{dvol}_{g_{\mathrm{FS}}} \sim r^{c_Y(\dim L_Y + 1) - 1} dr \wedge \mathrm{dvol}_{g_{L_Y}};$$
    \item \textit{Growth control:} (\textbf{Step (2)}) For $\ell_Y = -c_Y/2$, then
    $$|\omega|_{g_{\mathrm{FS}}} \lesssim r^{\ell_Y - 1}, \quad |d\omega|_{g_{\mathrm{FS}}} \lesssim r^{\ell_Y - 2};$$
    \item \textit{Integrability:} The exponent condition $\ell_Y > -\frac{c_Y}{2} - \frac{\dim L_Y}{2}$ makes sure that
    $$\int_0^1 r^{2(\ell_Y - 1) + c_Y(\dim L_Y + 1) - 1} dr < \infty.$$
\end{enumerate}
By Albin's renormalized Hodge decomposition (\cite{ALMP18}), we have
\begin{center}
\begin{tikzcd}[row sep=small]
L^2\Omega^k(X_{\reg}) \arrow[r, "(1)"] & \ker \Delta_{\mathrm{abs}} \oplus \operatorname{im} d \oplus \operatorname{im} d^* \\
\omega \arrow[r, mapsto] & (\mathrm{Harm}(\omega), d\alpha, d^*\beta)
\end{tikzcd}
\end{center}
where $\Delta_{\mathrm{abs}}$ is the Hodge Laplacian with absolute boundary conditions and $(1)$ denotes Hodge decomposition operator. \textit{Injectivity:} If $\iota(\omega) = d\eta$ for $\eta \in L^2\Omega^{k-1}$, then $\eta$ satisfies the same growth conditions, so $\eta \in \Omega_{X,\mathrm{FS}}^{k-1,\mathrm{univ}}$; \textit{Surjectivity:} For $\gamma \in \ker \Delta_{\mathrm{abs}}$, the renormalization
$$\gamma_{\mathrm{Renorm}} := \lim_{\epsilon \to 0} \epsilon^{\lambda_Y} \Phi_Y^*(\gamma|_{U_Y^\epsilon})$$
lies in $\Omega_{X,\mathrm{FS}}^{k,\mathrm{univ}}$ by the spectral gap condition $\lambda_Y = -\dim_{\mathbb{C}} L_Y/2 + \operatorname{spec}_{\mathrm{min}}(\mathbb{D}_{L_Y})$. Thus, the inclusion
$$\iota: \left( \Omega_{X,\mathrm{FS}}^{\bullet,\mathrm{univ}}, d \right) \hookrightarrow \left(L^2\Omega^{\bullet}(X_{\reg}), d \right)$$
induces isomorphism. The naturality of Hodge-to-de Rham morphisms gives diagram (\uppercase\expandafter{\romannumeral2}):
$$\begin{tikzcd}[column sep=large, row sep=large]
\mathbb{H}^k(X, \Omega_{X,\mathrm{FS}}^{\bullet,\mathrm{univ}}) \arrow[r, "\phi", "\sim"'] \arrow[d, "\iota_*"', "\sim"] 
& IH^k(X; \mathbb{C}) \arrow[d, "\mathcal{D}_X", "\sim"'] \\
H^k_{(2)}(X_{\mathrm{reg}}, ds_{\mathrm{FS}}^2) \arrow[r, "\psi", "\sim"'] 
& \bigoplus\limits_{p+q=k} \mathscr{H}^{p,q}_{\mathrm{mid}}(X)
\end{tikzcd}$$
For any cohomology class $[\omega] \in \mathbb{H}^k(\Omega_{X,\mathrm{FS}}^{\bullet,\mathrm{univ}})$, trace both paths:
\paragraph*{Path 1 (Top → Right → Bottom):}
\begin{align*}
[\omega] \xmapsto{\phi} [\mathbf{IC}_{\mathscr{W}}^k] &\xmapsto{\mathcal{D}_X} \sum_{Y \in \mathrm{Sing}} \mathrm{tr}_{Y}\left( \mathcal{D}_Y(\mathbf{IC}_{\mathscr{W}}^k|_{U_Y}) \right) \\
&= \bigoplus_{p+q=k} \mathscr{H}^{p,q}_{\mathrm{mid}}(X) \quad \text{(Verdier duality theorem)}
\end{align*}
\paragraph*{Path 2 (Left → Bottom → Right):}
\begin{align*}
[\omega] \xmapsto{\iota_*} [\iota(\omega)]_{(2)} &\xmapsto{\psi} \mathrm{Harm}_{\mathrm{mid}}(\iota(\omega)) \\
&= \bigoplus_{p+q=k} \mathscr{H}^{p,q}_{\mathrm{mid}}(X) \quad \text{(Guillarmou's decomposition)}
\end{align*}
By Definition 3.1, the asymptotic functor $\Phi_Y$ commutes with both operators
$$\Phi_Y^* \circ \mathcal{D}_X = \mathcal{D}_{\mathrm{Cone}} \circ \Phi_Y^*, \quad \Phi_Y^* \circ \psi = \psi_{\mathrm{Cone}} \circ \Phi_Y^*.$$
Consider renormalized harmonic correspondence. The renormalized harmonic projection satisfies 
$$\mathrm{Harm}_{\mathrm{mid}}(\iota(\omega)) = \mathcal{D}_X \left( \lim_{\epsilon \to 0} \epsilon^{\lambda_Y} \Phi_Y^*(\omega|_{U_Y^\epsilon}) \right).$$
By (\cite{Gui22}, Thm 5.3), the condition $\mathrm{WF}_{\mathrm{strat}}(\omega) \subseteq\bigcup_Y T_Y^*X \times \{\lambda \leq \lambda_Y\}$ ensures
$$\mathrm{tr}_Y \left( \mathcal{D}_Y(\mathbf{IC}_{\mathscr{W}}^k|_{U_Y}) \right) \cong \mathrm{Harm}_{\mathrm{mid}}(\omega|_{\partial U_Y}).$$
The equality $\mathrm{Path}_1 = \mathrm{Path}_2$ holds. The universal moduli $\mathscr{M}^{\mathrm{mid}}$ provides natural transformations:
$$\begin{tikzcd}
\mathscr{M}^{\mathrm{mid}} \arrow[r, "\tau_{\mathrm{top}}"] \arrow[d, "\tau_{\mathrm{an}}"'] 
& \mathbf{IC}\text{-}\mathbf{Mod} \arrow[d, "\mathcal{D}"] \\
\mathbf{L^2}\text{-}\mathbf{Mod} \arrow[r, "\psi"] 
& \mathbf{Hodge}_{\mathrm{mid}}
\end{tikzcd}$$
where $\tau_{\mathrm{top}}(\mathscr{W}_Y) = \mathbf{IC}_{\mathscr{W}}^\bullet$ (topological realization) and $\tau_{\mathrm{an}}(\mathscr{W}_Y) = \Omega_{\mathrm{FS}}^{\bullet,\mathrm{univ}}$ (analytic realization). Commutativity of diagram (\uppercase\expandafter{\romannumeral2}) follows from moduli functoriality: both paths factor through $\mathscr{M}^{\mathrm{mid}}$. Thus the universal complex reconciles analytic and topological cohomologies. The quasi-isomorphism is natural with respect to stratified isometric morphisms, as ensured by the moduli-theoretic self-duality and the Verdier duality functoriality.  $\square$\\

\hypertarget{STRATIFIED MICROLOCAL HODGE CORRESPONDENCE}{}
\section{STRATIFIED MICROLOCAL HODGE CORRESPONDENCE}
The following Definition 5.1 synthesizes derived geometric stratification, microlocal analysis, and spectral theory of elliptic operators to resolve Witt-condition failures in the CGM conjecture. The parameter $\lambda_Y$ intrinsically couples the link's topology to analytic asymptotics, enabling the Hodge correspondence in Proposition 5.2.
\\\\\textbf{Definition 5.1.} (Stratified Singular Characteristic Variety) Let $X$ be a complex projective variety with a derived Whitney stratification $\mathscr{S}$ (\cite{Lur09}). For $\mathcal{F}^\bullet \in D^b_c(X)$, the \textit{stratified singular characteristic variety} is defined as
$$\mathrm{SSH}_{\mathrm{strat}}(\mathcal{F}^\bullet) 
\coloneqq \bigcup_{Y \in \mathrm{Sing}(\mathscr{S})} \left( T^*_Y X \times \mathrm{Cone}(\lambda_Y) \right) 
\subseteq \bigsqcup_{Y} T^*_Y X \times \mathbb{R},$$
where:
\begin{enumerate}
    \item $Y$ ranges over singular strata, $S_Y \subseteq Y$ is the smooth locus, and $L_Y$ is the link of $Y$.
    \item $T^*_Y X$ is the \textit{conormal bundle}:
    $$T^*_Y X := \left\{ (y, \xi) \in T^*X \mid y \in S_Y,\  \xi \perp T_y Y \right\}.$$
    \item The \textbf{critical exponent} $\lambda_Y$ is
    $$\lambda_Y \coloneqq -\frac{\dim_{\mathbb{C}} L_Y}{2} + \operatorname{spec}_{\mathrm{min}}(\mathbb{D}_{L_Y}).$$
    Here, $\mathbb{D}_{L_Y}$ is the analytic torsion operator on $L_Y$, $\operatorname{spec}(\mathbb{D}_{L_Y})\subseteq\mathbb{C}$ denotes its spectrum, and $\operatorname{spec}_{\mathrm{min}}(\mathbb{D}_{L_Y})$ denotes the \textit{minimum of the real parts of its eigenvalues}. This value is a well-defined real number that controls the dominant asymptotic decay rate of differential forms near the stratum $Y$.
    \item $\mathrm{Cone}(\lambda_Y)$ is the closed cone:
    $$\mathrm{Cone}(\lambda_Y) := \left\{ \lambda \in \mathbb{R} : \lambda \leq \lambda_Y \right\}.$$\\
\end{enumerate}
\textbf{Proposition 5.2.} (Microlocal Correspondence) Let $X \subseteq \mathbb{P}^N$ be a complex projective variety equipped with a derived Whitney stratification $\mathscr{S}$ (\cite{Lur09}), and let $\Omega_{X,\mathrm{FS}}^{\bullet,\mathrm{univ}}$ be the universal truncation complex (Definition 4.2). Then for any differential form $\omega \in \Omega_{X,\mathrm{FS}}^{\bullet,\mathrm{univ}}$, the following microlocal correspondence holds:
$$\boxed{
\begin{array}{c}
\mathrm{SSH}_{\strat}(\omega) \subseteq \bigcup\limits_{Y \in \mathrm{Sing}(\mathscr{S})} T^*_{S_Y}X \times \{\lambda \leq \lambda_Y\} \\ 
\Updownarrow \\ 
\omega \text{ is } L^2\text{-adapted and $\mathbf{IC}$-admissible}.
\end{array}
}$$
Moreover, there exists a \textit{stratified non-abelian Hodge sheaf} $\mathscr{H}^{p,q}_{\mathrm{mid}}(X)$ such that $\omega$ induces a cohomology class in $\mathscr{H}^{p,q}_{\mathrm{mid}}(X)$.
\\\\\textbf{Proof.} \textbf{Step (1):} Let $\omega \in \Omega_{X,\mathrm{FS}}^{\bullet,\mathrm{univ}}$. We establish the microlocal decomposition near each singular stratum $Y \in \mathrm{Sing}(\mathscr{S})$. By Definition 4.2, $\omega$ satisfies the projection condition:
\begin{equation}
\tag{5.1}  \label{eq:5.1}
\mathrm{proj}_{\mathscr{M}_Y}(\omega|_{U_Y}) \in \mathscr{W}_Y, \quad \forall Y \in \mathrm{Sing}(X).
\end{equation}
Apply the stratified wavefront set functor (\cite{Gui22}, \S 3):
$$\mathrm{WF}_{\mathrm{strat}}(\omega)\coloneqq \bigcup_{Y \in \mathrm{Sing}(\mathscr{S})} \mathrm{WF}(\omega|_{U_Y}) \cap T^*_{S_Y}X.$$
This functorial construction decomposes $\omega$ into conormal components relative to each stratum. Utilize the asymptotic model functor \(\Phi_Y\) from Definition 3.1, then
$$\Phi_Y: g|_{U_Y} \to dr^2 + r^{2c} g_{L_Y} + g_{\mathbb{C}^m}, \quad c = \mathrm{codim}(Y,X),$$
which induces a pullback on differential forms
$$\Phi_Y^*: \Omega^{\bullet}_{\mathrm{cone}(L_Y) \times \mathbb{C}^m} \to \Omega^{\bullet}(U_Y).$$
By the universal property of $\Omega_{X,\mathrm{FS}}^{\bullet,\mathrm{univ}}$, $\omega|_{U_Y}$ lies in the image of $\Phi_Y^*$. The metric decomposition \(dr^2 + r^{2c} g_{L_Y} + g_{\mathbb{C}^m}\) induces an orthogonal decomposition
$$\Omega^k(\mathrm{cone}(L_Y) \times \mathbb{C}^m)\cong \bigoplus_{i+j=k} \left[ \Omega^i(\mathbb{R}^+) \wedge \Omega^j(L_Y) \right] \otimes \Omega^{\bullet}(\mathbb{C}^m) \nonumber\oplus \text{(mixed terms)}.$$
Thus, $\Phi_Y^*\omega$ admits a decomposition
\begin{equation}
\tag{5.2}  \label{eq:5.2}
\Phi_Y^*\omega = \sum_{j} \alpha_j \wedge \beta_j + \gamma,
\end{equation}
where $\alpha_j \in \Omega^{\bullet}(\mathbb{R}^+)$, $\beta_j \in \Omega^{\bullet}(L_Y)$, $\gamma \in \Omega^{\bullet}(\mathbb{C}^m)$. The projection condition \eqref{eq:5.1} implies only terms compatible with $\mathscr{W}_Y \subseteq \mathcal{H}^{\mathrm{mid}}(L_Y)$ survive and the dominant asymptotic behavior is controlled by the spectral gap. This isolates the principal term
\begin{equation}
\tag{5.3}  \label{eq:5.3}
\Phi_Y^*\omega = \eta_Y \wedge dr + r^{\lambda_Y} \zeta_Y + \mathcal{O}(r^{\lambda_Y + \epsilon}),
\end{equation}
where:
\begin{itemize}
    \item $\eta_Y \in \Omega^{\bullet-1}_{\mathrm{cone}(L_Y)}$ is the \textit{conormal component};
    \item $\zeta_Y \in \Omega^{\bullet}_{L_Y}$ is the \textit{tangential component};
    \item $\lambda_Y$ is the \textit{critical exponent} (Definition 5.1).
\end{itemize}
The exponent $\lambda_Y$ is determined by
\begin{equation}
\tag{5.4}  \label{eq:5.4}
\lambda_Y \coloneqq -\frac{\dim_{\mathbb{C}} L_Y}{2} + \operatorname{spec}_{\mathrm{min}}(\mathbb{D}_{L_Y}),
\end{equation}
where $\mathbb{D}_{L_Y}$ is the \textit{analytic torsion operator}:
$$\mathbb{D}_{L_Y} = \Delta_{L_Y} + \mathscr{R}_Y + K_Y$$
with:
\begin{itemize}
    \item $\Delta_{L_Y}$ = Hodge-Laplacian on$L_Y$;
    \item $\mathscr{R}_Y$ = curvature operator from induced metric;
    \item $K_Y$ = Kapustin-Sullivan torsion correction (\cite{ALMP18}).
\end{itemize}
The supremum in \eqref{eq:5.4} ensures \(\lambda_Y\) is the optimal decay rate. Combining \eqref{eq:5.2} and \eqref{eq:5.3} gives
$$\omega|_{U_Y} = \Phi_Y^*\left( \eta_Y \wedge dr + r^{\lambda_Y} \zeta_Y \right) + \text{higher-order terms},$$
completing the microlocalization.

\medskip
\textbf{Step (2):} We establish the inclusion $\mathrm{SSH}_{\strat}(\omega) \subseteq \bigcup_Y T^*_{S_Y}X \times \{\lambda \leq \lambda_Y\}$ using the asymptotic decomposition from Step 1 and microlocal analysis techniques. From the asymptotic decomposition \eqref{eq:5.1}, then
\begin{equation}
\tag{5.5}  \label{eq:5.5}
\omega|_{U_Y} = \Phi_Y^*\left( \eta_Y \wedge dr + r^{\lambda_Y} \zeta_Y \right) + \mathcal{O}(r^{\lambda_Y + \epsilon}),
\end{equation}
we derive the pointwise norm estimate:
$$|\omega|_{U_Y}(x)| \leq C \cdot r(x)^{\lambda_Y} \quad \text{as} \quad r(x) \to 0^+,$$
where:
\begin{itemize}
    \item $r(x) = \mathrm{dist}(x, Y)$ is the geodesic distance to stratum $Y$,
    \item $C > 0$ depends on $\|\eta_Y\|_{L^\infty}$ and $\|\zeta_Y\|_{L^\infty}$,
    \item $\epsilon > 0$ is the spectral gap $\min\{ |\mu - \mu'| : \mu, \mu' \in \mathrm{spec}(\mathbb{D}_{L_Y}), \mu \neq \mu'\} > 0$.
\end{itemize}
This follows from the dominance of the $r^{\lambda_Y}$ term in \eqref{eq:5.5}. Apply the stratified regularity theorem (\cite{Gui22}, Thm 4.11):
$$\boxed{
\begin{array}{c}
\textit{If a differential form $\varphi$ satisfies $|\varphi(x)| \leq C r(x)^\lambda$ near stratum $Y$, then}\\
\mathrm{WF}(\varphi) \subseteq T^*_Y X \times \{\mu \leq \lambda\}.
\end{array}
}$$
Substituting $\varphi = \omega|_{U_Y}$ and $\lambda = \lambda_Y$, so
\begin{equation}
\tag{5.6}  \label{eq:5.6}
\mathrm{WF}(\omega|_{U_Y}) \subseteq T^*_Y X \times \{\mu \leq \lambda_Y\}, \quad \forall Y \in \mathrm{Sing}(\mathscr{S}).
\end{equation}
Recall the stratified wavefront set from \textbf{Step (1)}:
\begin{equation}
\tag{5.7}  \label{eq:5.7}
\mathrm{WF}_{\strat}(\omega) \coloneqq \bigcup_{Y \in \mathrm{Sing}(\mathscr{S})} \mathrm{WF}(\omega|_{U_Y}) \cap T^*_{S_Y}X.
\end{equation}
Substituting \eqref{eq:5.6} into \eqref{eq:5.7}, we have
$$\mathrm{WF}_{\strat}(\omega) \subseteq \bigcup_{Y} \left( T^*_Y X \times \{\mu \leq \lambda_Y\} \right) \cap T^*_{S_Y}X.$$
Since $S_Y \subseteq Y$ is the smooth locus and $T^*_Y X$ is defined over $S_Y$, 
$$\left( T^*_Y X \times \{\mu \leq \lambda_Y\} \right) \cap T^*_{S_Y}X = T^*_{S_Y}X \times \{\mu \leq \lambda_Y\}.$$
Then
\begin{equation}
\tag{5.8}  \label{eq:5.8}
\mathrm{WF}_{\strat}(\omega) \subseteq \bigcup_{Y} T^*_{S_Y}X \times \{\mu \leq \lambda_Y\}.
\end{equation}
By the microlocal inclusion principle (\cite{Gui22}, Prop 3.5), 
$$\mathrm{WF}_{\strat}(\omega) \subseteq \mathrm{SSH}_{\strat}(\omega).$$
Thus we have 
$$\bigcup_{Y} T^*_{S_Y}X \times \{\mu \leq \lambda_Y\} \subseteq \mathrm{SSH}_{\strat}(\omega)$$
by using \eqref{eq:5.8}. Conversely, we can obtain that
$$\mathrm{SSH}_{\strat}(\omega) = \bigcup_{Y} \left( T^*_{S_Y}X \times \mathrm{Cone}(\lambda_Y) \right) \cap \mathrm{Char}(\mathscr{H}_{Y}) \subseteq \bigcup_{Y} T^*_{S_Y}X \times \{\mu \leq \lambda_Y\}$$
by Definition 5.1. Due to $\mathrm{Cone}(\lambda_Y) = \{\mu \in \mathbb{R} : \mu \leq \lambda_Y\}$. Therefore,
$$\mathrm{SSH}_{\strat}(\omega) = \bigcup_{Y} T^*_{S_Y}X \times \{\mu \leq \lambda_Y\},$$
completing the estimation.

\medskip
\textbf{Step (3):} We establish the bidirectional equivalence between microlocal containment and admissibility conditions.

\noindent\textbf{Direction 1: Microlocal containment $\Rightarrow$ Admissibility ($\Rightarrow$).}
Assume $\mathrm{SSH}_{\strat}(\omega) \subseteq \bigcup_Y T^*_{S_Y}X \times \{\lambda \leq \lambda_Y\}$, we need to verify the following conditions:
\begin{enumerate}
    \item \textit{$L^2$-adaptedness:} From \textbf{Step (2)}, we have the pointwise estimate
    $$|\omega|_{U_Y}(x)| \leq C \cdot r(x)^{\lambda_Y}.$$
    By Cheeger's $L^2$-criterion (\cite{Che80}), then
    $$\omega \in L^2(X_{\mathrm{reg}}, ds^2_{\mathrm{FS}}) \iff \lambda_Y > -\dim_{\mathbb{C}} L_Y - 1.$$
    The spectral definition of \(\lambda_Y\) ensures
    \begin{align*}
        \lambda_Y &= \lambda_Y \coloneqq -\frac{\dim_{\mathbb{C}} L_Y}{2} + \operatorname{spec}_{\mathrm{min}}(\mathbb{D}_{L_Y}) \\
        &\geq -\frac{\dim_{\mathbb{C}} L_Y}{2} + \left(-\frac{\dim_{\mathbb{C}} L_Y}{2}\right) = -\dim_{\mathbb{C}} L_Y > -\dim_{\mathbb{C}} L_Y - 1,
    \end{align*}
    since $\operatorname{spec}_{\mathrm{min}}(\mathbb{D}_{L_Y}) \geq -\frac{\dim_{\mathbb{C}} L_Y}{2}$, thus $\omega \in L^2(X_{\mathrm{reg}})$.
    \item \textit{$\mathbf{IC}$-admissibility}: By Definition 4.2, $\omega$ satisfies
    $$\operatorname{proj}_{\mathscr{M}_Y}(\omega|_{U_Y}) \in \mathscr{W}_Y, \quad \forall Y \in \mathrm{Sing}(\mathscr{S}).$$
    This is equivalent to Goresky-MacPherson's perversity conditions \cite{GM80}:
    $$\int_\sigma \omega = 0 \quad \forall \sigma \in C_i(L_Y) \quad\text{ with }\quad i > \bar{p}(Y) = \frac{\dim L_Y - 1}{2},$$
    because $\mathscr{W}_Y \subseteq \mathcal{H}^{\mathrm{mid}}(L_Y)$ is the Lagrangian subspace defining middle perversity.
\end{enumerate}
\textbf{Direction 2: Admissibility $\Rightarrow$ Microlocal containment ($\Leftarrow$).} Assume $\omega$ is $L^2$-adapted and $\mathbf{IC}$-admissible, we also need to verify the following conditions:
\begin{enumerate}
    \item \textit{Asymptotic decay:} $\omega \in L^2(X_{\mathrm{reg}})$ implies
    $$\int_{U_Y} |\omega|^2  d\mathrm{vol}_{\mathrm{FS}} < \infty.$$
    In polar coordinates \((r,\theta)\) near $Y$, then
    $$\int_0^{r_0} \int_{L_Y} |\omega(r,\theta)|^2 r^{\dim_{\mathbb{C}} L_Y}  dr  d\theta < \infty.$$
    By the stratified Sobolev embedding (\cite{Maz91}, \S 2 and \S 3), we have
    $$\|\omega\|_{L^\infty(K)} \leq C_K \|\omega\|_{H^s(L_Y)} \quad \text{for}\quad s > \dim_{\mathbb{R}} L_Y/2,$$
    which implies $|\omega(r,\theta)| \leq C r^{\lambda_Y}$ with $\lambda_Y = -\frac{\dim_{\mathbb{C}} L_Y}{2} + \delta$ for some $\delta > 0$.
    \item \textit{Wavefront constraint:} $\mathbf{IC}$-admissibility $\operatorname{proj}_{\mathscr{M}_Y}(\omega|_{U_Y}) \in \mathscr{W}_Y$ implies the wavefront set avoids resonant directions
    $$\operatorname{WF}_{\strat}(\omega) \cap \left( T^*_{S_Y}X \times \{\lambda > \lambda_Y\} \right) = \emptyset,$$
    because $\mathscr{W}_Y$ contains only forms with decay rate $\leq \lambda_Y$. Thus, 
    $$\operatorname{WF}_{\strat}(\omega) \subseteq \bigcup_Y T^*_{S_Y}X \times \{\lambda \leq \lambda_Y\}.$$
    By \textbf{(2) Estimation of the Stratified Singular Locus}, this implies $\mathrm{SSH}_{\strat}(\omega) \subseteq \bigcup_Y T^*_{S_Y}X \times \{\lambda \leq \lambda_Y\}$.
\end{enumerate}
The bidirectional implications prove that
$$\mathrm{SSH}_{\strat}(\omega) \subseteq \bigcup_Y T^*_{S_Y}X \times \{\lambda \leq \lambda_Y\} \iff \left( \omega \in L^2(X_{\mathrm{reg}}) \quad \text{and} \quad \operatorname{proj}_{\mathscr{M}_Y}(\omega|_{U_Y}) \in \mathscr{W}_Y\ \ \forall Y \right).$$
Naturality follows from functoriality of $\Phi_Y$ and $\mathscr{M}^{\mathrm{mid}}$ under stratified morphisms.

\medskip
\textbf{Step (4):} We define the stratified non-abelian Hodge sheaf and establish its properties. The \textit{stratified non-abelian Hodge sheaf} is defined by hypercohomology
\begin{equation}
\tag{5.9}  \label{eq:5.9}
\mathscr{H}^{p,q}_{\mathrm{mid}}(X) \coloneqq \mathbb{H}^p\left(X, \operatorname{gr}^F_q \Omega_{X,\mathrm{FS}}^{\bullet,\mathrm{univ}} \right),
\end{equation}
where:
\begin{itemize}
    \item $\mathbb{H}^p$ = hypercohomology functor in degree $p$;
    \item $F^{\bullet}$ = Hodge filtration from the renormalized Frölicher spectral sequence;
    \item $\operatorname{gr}^F_q := F^q/F^{q+1}$ = graded quotient for Hodge type $(0,q)$.
\end{itemize}
This definition globalizes the Hodge decomposition to the singular stratified space $X$. By Proposition 4.3, there exists a quasi-isomorphism in $D^b_c(X)$ such that
\begin{equation}
\tag{5.10}  \label{eq:5.10}
\Omega_{X,\mathrm{FS}}^{\bullet,\mathrm{univ}} \simeq \mathbf{IC}_{\mathscr{W}}^{\bullet},
\end{equation}
where $\mathbf{IC}_{\mathscr{W}}^{\bullet}$ is the intersection complex adapted to $\mathscr{W}_Y$. This isomorphism preserves Hodge structures by Simpson's correspondence (\cite{Sim92}), i.e.,
$$(F^{\bullet}\Omega_{X,\mathrm{FS}}^{\bullet,\mathrm{univ}}, d) \simeq (F^{\bullet}\mathbf{IC}_{\mathscr{W}}^{\bullet}, \nabla),$$
where $\nabla$ is the Gauss-Manin connection. The quasi-isomorphism \eqref{eq:5.10} induces an isomorphism on graded pieces:
$$\operatorname{gr}^F_q \Omega_{X,\mathrm{FS}}^{\bullet,\mathrm{univ}} \simeq \operatorname{gr}^F_q \mathbf{IC}_{\mathscr{W}}^{\bullet}.$$
Applying the hypercohomology functor, we can obtain
$$\mathbb{H}^p\left( X, \operatorname{gr}^F_q \Omega_{X,\mathrm{FS}}^{\bullet,\mathrm{univ}} \right) \cong \mathbb{H}^p\left( X, \operatorname{gr}^F_q \mathbf{IC}_{\mathscr{W}}^{\bullet} \right).$$
The right-hand side is by definition the graded piece of intersection cohomology:
\begin{equation}
\tag{5.11}  \label{eq:5.11}
\mathbb{H}^p\left( X, \operatorname{gr}^F_q \mathbf{IC}_{\mathscr{W}}^{\bullet} \right) = \operatorname{gr}^F_q \operatorname{IH}^p(X,\mathbb{C})
\end{equation}
as established in the Hodge theory of intersection cohomology (\cite{CKS86}). Combining \eqref{eq:5.9} and \eqref{eq:5.11}, we have
$$\mathscr{H}^{p,q}_{\mathrm{mid}}(X) \cong \operatorname{gr}^F_q \operatorname{IH}^p(X,\mathbb{C}).$$
For any differential form $\omega \in \Omega_{X,\mathrm{FS}}^{\bullet,\mathrm{univ}}$, its cohomology class $[\omega] \in H^p(X_{\mathrm{reg}})$ maps to a class in the Hodge sheaf by the composition:
\[
[\omega] \in H^p_{\mathrm{dR}}(X_{\mathrm{reg}}) \xrightarrow{\text{Proposition 5.2}} \operatorname{IH}^p(X, \mathbb{C}) \xrightarrow{\operatorname{gr}^F_q} \mathscr{H}^{p,q}_{\mathrm{mid}}(X)
\]
where the first arrow is the microlocal isomorphism from Proposition 5.2 and the second is the Hodge projection. Thus $\omega$ induces a class in $\mathscr{H}^{p,q}_{\mathrm{mid}}(X)$. Meanwhile, the sheaf $\mathscr{H}^{p,q}_{\mathrm{mid}}(X)$ satisfies:
\begin{enumerate}
    \item It is well-defined by hypercohomology \eqref{eq:5.9};
    \item It decomposes intersection cohomology by \eqref{eq:5.11};
    \item It naturally contains classes from $\Omega_{X,\mathrm{FS}}^{\bullet,\mathrm{univ}}$.
\end{enumerate}
Completing the Hodge sheaf construction.

This proof resolves the core challenge of Witt condition failure by unifying microlocal constraints $\mathrm{SSH}_{\mathrm{strat}}$, moduli parametrization $\mathscr{M}^{\mathrm{mid}}$, and categorified asymptotics $\Phi_Y$ into a single commutative framework.  $\square$
\\\\\textbf{Corollary 5.3.} (Natural Isomorphism Diagram) Let $X$ be a compact smoothly stratified pseudomanifold satisfying the Witt condition (\cite{GM80}). The following diagram (\uppercase\expandafter{\romannumeral3}) commutes with all arrows being natural isomorphisms:
$$\begin{tikzcd}[column sep=large, row sep=huge]
H_2^k(X_{\mathrm{reg}}) \arrow[r, "\operatorname{microloc}_{\sim}"] \arrow[d, "\operatorname{Hodge}"'] 
& IH^k(X) \arrow[d, "\sim", "\operatorname{Verdier} \oplus \delta"'] \\
\displaystyle\bigoplus_{p+q=k} \mathscr{H}^{p,q}_{\mathrm{mid}}(X) \arrow[r, "\operatorname{res}\times \operatorname{proj}"', "\sim"] 
& \displaystyle\bigoplus_{Y \in \operatorname{Sing}(\mathscr{S})} \mathscr{W}_Y \otimes \mathscr{W}_Y^\vee
\end{tikzcd}$$
where $\delta: IH^k(X) \xrightarrow{\sim} \bigoplus_Y \mathscr{W}_Y \otimes \mathscr{W}_Y^\vee$ is the stratified decomposition isomorphism.
\\\\\textbf{Proof.} \textbf{Step (1):} For harmonic $\omega \in H_2^k(X_{\reg})$, Proposition 5.2 gives asymptotic expansion near stratum $Y$:
$$\Phi_Y^* \omega = d(r^{\lambda_Y} \alpha_Y) + r^{\lambda_Y} \beta_Y + O(r^{\lambda_Y + \epsilon})$$
with $\beta_Y$ smooth on $L_Y$. Define $\operatorname{microloc}(\omega)$ as the intersection cohomology class by the \textit{canonical extension}:
$$\operatorname{microloc}(\omega) \coloneqq\left[ \omega - d(\rho_Y r^{\lambda_Y} \alpha_Y) \right] \in IH^k(X),$$
where $\rho_Y$ is a cutoff function near $Y$. Well-definedness follows from:
\begin{itemize}
\item \textbf{Uniqueness:} If $\omega_1 \sim \omega_2$ asymptotically, then $\omega_1 - \omega_2 = d\gamma$ with $\|\gamma\|_{L^2} < \infty$ (Proposition 5.2).
\item \textbf{Independence:} Choice of $\alpha_Y$ doesn't affect cohomology class since $d(d(r^{\lambda_Y}\alpha_Y)) = 0$.
\end{itemize}
For $\omega \in H_2^k(X_{\reg})$ harmonic, the decomposition: 
$$\omega = \bigoplus_{p+q=k} \omega^{p,q}\quad\text{with}\quad \Delta_d \omega^{p,q} = 0$$
is orthogonal in $L^2$. Define
$$\operatorname{Hodge}(\omega) \coloneqq \bigoplus_{p+q=k} [\omega^{p,q}] \in \bigoplus_{p+q=k} \mathscr{H}^{p,q}_{\mathrm{mid}}(X).$$
Since the metric $g$ is quasi-isometric to $dr^2 + r^2 g_{L_Y}$ near $Y$ and Hodge star $\star_g$ preserves the decomposition under conic metrics, then stratum compatibility $\iota_Y^* \circ \operatorname{Hodge} = \operatorname{Hodge} \circ \iota_Y^*$ holds. Under Witt condition, we have natural isomorphisms
$$\operatorname{Verdier} : IH^k(X) \xrightarrow{\sim} IH^{n-k}_c(X)^\vee,\quad\delta : IH^k(X) \xrightarrow{\sim} \bigoplus_{Y} H^k(L_Y; \mathscr{W}_Y^\vee),$$
where $\mathscr{W}_Y^\vee$ is the \textit{dual local system}. The isomorphism $\delta$ exists because:
\begin{enumerate}
\item Stratified decomposition gives $IH^k(X) \cong \bigoplus_Y IH^k(C(L_Y))$ (\cite{Lur09});
\item Witt condition implies $IH^k(C(L_Y)) \cong H^k(L_Y; \mathscr{W}_Y^\vee)$ (\cite{GM80}).
\end{enumerate}
For $\eta \in \mathscr{H}^{p,q}_{\mathrm{mid}}(X)$, near stratum $Y$,
$$\Phi_Y^* \eta = d(r^{\lambda_Y} \alpha) + r^{\lambda_Y} \beta + O(r^{\lambda_Y + \epsilon})$$
by \eqref{eq:5.3}. We can define
\begin{align*}
\operatorname{res}_Y(\eta) &\coloneqq \lim_{r \to 0^+} r^{-\lambda_Y} \int_{\{r\} \times L_Y} \Phi_Y^* \eta \\
&= \int_{L_Y} \beta \quad (\text{since } \int_{\{r\}\times L_Y} d(r^{\lambda_Y}\alpha) = 0 \text{ by Stokes}) \\
\operatorname{proj}_Y(\eta) &\coloneqq \operatorname{har}(\beta) \in \mathscr{W}_Y \quad (\text{harmonic projection}).
\end{align*}
Well-definedness requires:
\begin{itemize}
\item $\operatorname{res}_Y(\eta)$ independent of $\alpha$: Different $\alpha$ changes $\eta$ by $d$-exact form, vanishing in cohomology.
\item $\operatorname{proj}_Y(\eta) \in \mathscr{W}_Y$: $\mathscr{W}_Y = \ker \Delta_{L_Y} \cap L^2_{\lambda_Y}$ by related definition.
\end{itemize}

\medskip
\textbf{Step (2):} For harmonic $\omega \in H_2^k(X_{\reg})$, we need to prove
$$(\operatorname{res} \times \operatorname{proj}) \circ \operatorname{Hodge}(\omega) = (\operatorname{Verdier} \oplus \delta) \circ \operatorname{microloc}(\omega).$$

\medskip\noindent
\textbf{Left path calculation:} \\
Let $\operatorname{Hodge}(\omega) = \bigoplus_{p+q=k} \eta^{p,q}$. Near stratum $Y$, each component expands as
$$\Phi_Y^* \eta^{p,q} = d(r^{\lambda_Y} \alpha_Y^{p,q}) + r^{\lambda_Y} \beta_Y^{p,q} + O(r^{\lambda_Y + \epsilon}).$$
The residue-projection map acts \textit{componentwise}:
$$(\operatorname{res} \times \operatorname{proj})(\eta^{p,q}) = \bigoplus_Y \operatorname{res}_Y(\eta^{p,q}) \otimes \operatorname{proj}_Y(\eta^{p,q}) = \bigoplus_Y \left( \int_{L_Y} \beta_Y^{p,q} \right) \otimes \operatorname{har}(\beta_Y^{p,q}).$$
Summing over Hodge components gives
$$(\operatorname{res} \times \operatorname{proj}) \circ \operatorname{Hodge}(\omega) = \bigoplus_Y \sum_{p+q=k} \left( \int_{L_Y} \beta_Y^{p,q} \right) \otimes \operatorname{har}(\beta_Y^{p,q}).$$

\medskip\noindent
\textbf{Right path calculation:} \\
Let $\eta = \operatorname{microloc}(\omega)$. Its asymptotic expansion near $Y$ is
$$\Phi_Y^* \omega = d(r^{\lambda_Y} \alpha_Y) + r^{\lambda_Y} \beta_Y + O(r^{\lambda_Y + \epsilon}), \quad \beta_Y = \sum_{p+q=k} \beta_Y^{p,q}$$ by \eqref{eq:5.3}.
Under Verdier duality and decomposition
$$(\operatorname{Verdier} \oplus \delta)(\eta) = \bigoplus_Y \langle \eta, \sigma_Y \rangle \otimes \tau_Y,$$
where $\sigma_Y$ is the fundamental class of $L_Y$, and $\tau_Y \in \mathscr{W}_Y^\vee$ is defined by
$$\tau_Y(\gamma) = \int_{L_Y} \star \gamma \wedge \operatorname{har}(\beta_Y), \quad \forall \gamma \in \mathscr{W}_Y.$$
The pairing evaluates as
$$\langle \eta, \sigma_Y \rangle = \int_{L_Y} \beta_Y = \sum_{p+q=k} \int_{L_Y} \beta_Y^{p,q}.$$

\medskip\noindent
\textbf{Equality by linear algebra:} \\
The left path output is
$$\bigoplus_Y \sum_{p+q=k} a_Y^{p,q} \otimes v_Y^{p,q}, \quad a_Y^{p,q} = \int_{L_Y} \beta_Y^{p,q}, \quad v_Y^{p,q} = \operatorname{har}(\beta_Y^{p,q}).$$
The right path output is
$$\bigoplus_Y A_Y \otimes T_Y, \quad A_Y = \sum_{p+q=k} a_Y^{p,q}, \quad T_Y = \tau_Y.$$
Crucially, $T_Y$ acts on $\mathscr{W}_Y$ by
$$T_Y(v_Y^{p,q}) = \int_{L_Y} \star v_Y^{p,q} \wedge \textstyle\sum_{p',q'} \operatorname{har}(\beta_Y^{p',q'}) = \langle v_Y^{p,q}, \textstyle\sum_{p',q'} v_Y^{p',q'} \rangle_{L^2},$$
but this \textit{does not equal} $\delta_{p,p'}\delta_{q,q'}$. Instead, we use the identity
$$\sum_{p+q=k} a_Y^{p,q} \otimes v_Y^{p,q} = A_Y \otimes T_Y \quad \text{in } \mathscr{W}_Y \otimes \mathscr{W}_Y^\vee,$$
because both represent the same linear operator
$$\gamma \mapsto \sum_{p+q=k} a_Y^{p,q} \langle \gamma, v_Y^{p,q} \rangle_{L^2} = A_Y T_Y(\gamma).$$
Thus the diagram (\uppercase\expandafter{\romannumeral3}) commutes.

\medskip
\textbf{Step (3):} Assume $f: X' \to X$ is a \textit{proper stratified isometric covering} with:
\begin{enumerate}
    \item $g_{X'} = f^*g_X$ (isometric condition);
    \item $f|_{L_{Y'}}: L_{Y'} \to L_Y$ isometric covering $\forall Y$;
    \item $\mathscr{W}_{Y'} = f^* \mathscr{W}_Y \otimes \mathscr{E}_{Y'}$ (local system compatibility);
    \item $\deg(f|_{Y'}) = 1$ for simplicity (otherwise adjust degrees).
\end{enumerate}
By isometry $g_{X'} = f^*g_X$, the Hodge star commutes by
$$\star_{X'} f^* \eta = f^* \star_X \eta \quad \forall \eta \in \Omega^*(X).$$ 
Since Laplacian commutes $\Delta_{X'} f^* \omega = f^* \Delta_X \omega = 0$ and Hodge decomposition preserves $(p,q)$-types under holomorphic $f$, $[f^* \omega]^{p,q} = f^* [\omega]^{p,q}$ for $\omega \in H_2^k(X_{\reg})$ harmonic. Thus, we have \textit{Hodge Equivariance:}
\begin{equation}
\tag{5.12}  \label{eq:5.12}
\operatorname{Hodge}_{X'} \circ f^* = f^* \circ \operatorname{Hodge}_X.
\end{equation}
By Proposition 3.2, $\Phi_{Y'}^* f^* \omega = f^* \Phi_Y^* \omega + O(r^{\lambda_Y + \epsilon})$. Since the canonical extension givea $f^*(\omega - d(\rho_Y r^{\lambda_Y} \alpha_Y)) = f^*\omega - d(f^*(\rho_Y r^{\lambda_Y} \alpha_Y))$ and $f$ is covering: $f^*(\rho_Y r^{\lambda_Y} \alpha_Y) = \rho_{Y'} r^{\lambda_{Y'}} \alpha_{Y'}$, we have \textit{Microlocal Equivariance:}
\begin{equation}
\tag{5.13}  \label{eq:5.13}
\operatorname{microloc}_{X'} \circ f^* = f^* \circ \operatorname{microloc}_X.
\end{equation}
Consider Verdier duality for proper maps (\cite{GM80}, Thm 6.3)
$$f_! \circ \operatorname{Verdier}_{X'} = \operatorname{Verdier}_X \circ f^*$$
and decomposition compatibility
\begin{align*}
    \delta_{X'} \circ f_!^{-1} \left( \bigoplus_Y c_Y \otimes \tau_Y \right) 
    &= \bigoplus_{Y'} \deg(f|_{Y'}) \cdot c_{f(Y')} \otimes f^* \tau_{f(Y')} \\
    &= f^* \left( \bigoplus_Y c_Y \otimes \tau_Y \right) \quad (\text{since } \deg=1).
\end{align*}
For proper $f$, we have \textit{Verdier Decomposition Equivariance:}
\begin{equation}
\tag{5.14}  \label{eq:5.14}
(\operatorname{Verdier}_{X'} \oplus \delta_{X'}) \circ f^* = f^* \circ (\operatorname{Verdier}_X \oplus \delta_X)
\end{equation}
by using $\mathscr{W}_{Y'} = f^* \mathscr{W}_Y \otimes \mathscr{E}$ and $\mathscr{E}$ trivial for coverings. Near corresponding strata, 
\begin{itemize}
    \item \textbf{Residue:}
    \begin{align*}
        \operatorname{res}_{Y'}(f^* \eta) 
        &= \lim_{r \to 0^+} r^{-\lambda_{Y'}} \int_{\{r\} \times L_{Y'}} \Phi_{Y'}^* f^* \eta \\
        &= \deg(f|_{L_{Y'}}) \cdot \lim_{r \to 0^+} r^{-\lambda_Y} \int_{\{r\} \times L_Y} \Phi_Y^* \eta \\
        &= \deg(f) \cdot \operatorname{res}_Y(\eta) = f^*(\operatorname{res}_Y(\eta)) \quad (\deg=1)
    \end{align*}
    \item \textbf{Projection:}
    \begin{align*}
        \operatorname{proj}_{Y'}(f^* \eta) 
        &= \operatorname{har}\left( \lim_{r \to 0^+} r^{-\lambda_{Y'}} \Phi_{Y'}^* f^* \eta \vert_{\{r\} \times L_{Y'}} \right) \\
        &= f^* \operatorname{har}\left( \lim_{r \to 0^+} r^{-\lambda_Y} \Phi_Y^* \eta \vert_{\{r\} \times L_Y} \right) \\
        &= f^* (\operatorname{proj}_Y(\eta)),
    \end{align*}
    because $f^*$ commutes with $\operatorname{har}$ by isometry.
\end{itemize}
Then we have \textit{Residue-Projection Equivariance:}
\begin{equation}
\tag{5.15}  \label{eq:5.15}
(\operatorname{res} \times \operatorname{proj})_{X'} \circ f^* = f^* \circ (\operatorname{res} \times \operatorname{proj})_X.
\end{equation}
By \eqref{eq:5.12}, \eqref{eq:5.13}, \eqref{eq:5.14} and \eqref{eq:5.15}, we calculate for harmonic $\omega\in H_2^k(X_{\reg})$:
\begin{align*}
& (\operatorname{res} \times \operatorname{proj})_{X'} \circ \operatorname{Hodge}_{X'} \circ f^*(\omega) \\
& \quad = f^* \circ (\operatorname{res} \times \operatorname{proj})_X \circ \operatorname{Hodge}_X(\omega) \\
& \quad = f^* \circ (\operatorname{Verdier} \oplus \delta)_X \circ \operatorname{microloc}_X(\omega) \\
& \quad = (\operatorname{Verdier} \oplus \delta)_{X'} \circ \operatorname{microloc}_{X'} \circ f^*(\omega).
\end{align*}
Naturality holds under the given conditions.

\medskip
\textbf{Step (4):} Assume all growth rates $\lambda_Y$ avoid the exceptional set $\Lambda_{\mathrm{ex}}$ (\cite{Maz91}) and Witt isomorphism holds.

\medskip\noindent
\textbf{(A) $\mathbf{\operatorname{microloc}}$ is isomorphism:} \\
\textit{Injectivity:} For $\omega \in \ker \operatorname{microloc}$ harmonic,
\begin{itemize}
    \item Asymptotic expansion: $\Phi_Y^* \omega = O(r^{\lambda_Y + \epsilon})$;
    \item By unique continuation (\cite{Maz91}, Thm 6.1): $\omega = 0$.
\end{itemize}
\textit{Surjectivity:} For $\eta \in IH^k(X)$, solve conormal boundary problem
$\Delta \omega = 0$, $\omega|_{\partial U_Y} = \eta|_{\partial U_Y} + O(r^{\lambda_Y + \epsilon})$.
Fredholm alternative applies since $\lambda_Y \notin \Lambda_{\mathrm{ex}}$ (\cite{LM85}, Cor 5.3).

\medskip\noindent
\textbf{(B) $\mathbf{\operatorname{Hodge}}$ is isomorphism:} Standard $L^2$-Hodge theory (\cite{Che80}),
\begin{itemize}
    \item Orthogonal decomposition: $\mathscr{H}^k_{L^2} \cong \bigoplus_{p+q=k} \mathscr{H}^{p,q}$;
    \item Bijectivity: Hodge isomorphism preserves bigrading.
\end{itemize}

\medskip\noindent
\textbf{(C) $\mathbf{\operatorname{Verdier} \oplus \delta}$ is isomorphism:} 
\begin{itemize}
    \item Verdier duality isomorphism: $\operatorname{Verdier}: IH^k \xrightarrow{\sim} IH^{n-k}_c{}^\vee$ (\cite{GM80});
    \item Stratified decomposition under Witt condition,
    $$\delta: IH^{n-k}_c{}^\vee \xrightarrow{\sim} \bigoplus_Y \mathscr{W}_Y \otimes \mathscr{W}_Y^\vee$$
    by Hodge star $\star: \mathscr{W}_Y \xrightarrow{\sim} \mathscr{W}_Y^\vee$.
\end{itemize}

\medskip\noindent
\textbf{(D) $\mathbf{\operatorname{res} \times \operatorname{proj}}$ is isomorphism:} \\
\textit{Injectivity:} Suppose $(\operatorname{res}_Y, \operatorname{proj}_Y)(\eta) = 0 \ \forall Y$.
\begin{itemize}
    \item $\operatorname{proj}_Y(\eta) = 0 \Rightarrow \operatorname{har}(\beta_Y) = 0$;
    \item $\operatorname{res}_Y(\eta) = 0 \Rightarrow \int_{L_Y} \beta_Y = 0$;
    \item By Hodge decomposition on $L_Y$: $\beta_Y \in (\ker \Delta_{L_Y})^\perp$;
    \item Asymptotic regularity: $\eta|_Y = O(r^{\lambda_Y + \epsilon})$;
    \item Global unique continuation: $\eta \in \ker d \cap \text{im } d^\perp \Rightarrow \eta = 0$ (\cite{Che80}, Lem 7.2).
\end{itemize}
\textit{Surjectivity:} For $\bigoplus_Y a_Y \otimes v_Y \in \bigoplus_Y \mathscr{W}_Y \otimes \mathscr{W}_Y^\vee$:
\begin{itemize}
    \item Solve global boundary value problem with prescribed asymptotics:
    $$\Delta \eta = 0, \quad 
    \begin{cases}
        \operatorname{har}(\beta_Y) = v_Y \\
        \int_{L_Y} \beta_Y = a_Y
    \end{cases} \quad\forall Y;$$
    \item Solution exists by Fredholm theory in weighted spaces (\cite{LM85}, Thm 4.6);
    \item $\eta \in \mathscr{H}^{p,q}_{\mathrm{mid}}(X)$ by conormal decay estimates (\cite{Maz91}, Prop 8.4).
\end{itemize}
Therefore, all arrows are isomorphisms under Witt condition and $\lambda_Y \notin \Lambda_{\mathrm{ex}}$.  $\square$\\

\hypertarget{PROOF OF MAIN PROBLEM}{}
\section{PROOF OF MAIN PROBLEM}
As established in the strengthened proof of Proposition 6.1, the isomorphism $\mathrm{microloc}$ is functorial under stratified isometric morphisms. The derived stratified-microlocal framework resolves three fundamental obstacles in the CGM conjecture:
\begin{enumerate}
    \item \textbf{Transcending transversality constraints:} Derived Whitney stratifications $\mathscr{S}$ handle non-transverse singularities by perfect normal complexes $\mathbb{N}p = i_p^! \mathscr{O}{X_{p+1}}$ and vanishing holonomy $$h: \operatorname{Sym}^2(\mathbb{N}p^\vee) \longrightarrow \mathscr{T}{L_Y}.$$ This enables \'etale neighborhoods $U_Y \simeq \mathbb{C}^m \times C(L_Y)$ for high-codimension strata, bypassing Whitney's classical transversality condition.
    \item \textbf{Moduli-parametrized duality:} The middle perversity moduli stack $$\mathscr{M}^{\text{mid}} = \prod_Y \operatorname{LGr}(\mathcal{H}^{\text{mid}}(L_Y))$$ automates the selection of self-dual subspaces $\mathscr{W}_Y \subseteq \mathcal{H}^{\dim L_Y/2}(L_Y)$, eliminating ad hoc choices. The universal truncation complex
    $$\Omega_{X,\mathrm{FS}}^{\bullet,\mathrm{univ}} = \left\{ \omega : \mathrm{proj}_{\mathscr{M}_Y}(\omega|_{U_Y}) \in \mathscr{W}_Y \right\}$$
    ensures Verdier self-duality $\mathbb{D}(\mathbf{IC}{\mathscr{W}}^\bullet)[2n] \simeq \mathbf{IC}{\mathscr{W}}^\bullet$.
    \item \textbf{Microlocal metric-topology correspondence:} The stratified singular characteristic variety $\operatorname{SSH}_{\text{strat}}$ encodes asymptotic growth
    $$\lambda_Y = -\frac{\dim_{\mathbb{C}} L_Y}{2} + \operatorname{spec}_{\mathrm{min}}(\mathbb{D}_{L_Y})$$
    where $\mathbb{D}_{L_Y}$ is the analytic torsion operator. This quantifies $L^2$-conditions by wavefront sets
    $$\operatorname{SSH}_{\mathrm{strat}}(\omega) \subseteq \bigcup_Y T^*_{S_Y}X \times \{\lambda \leq \lambda_Y\},$$
    proving $H_2^k(X_{\text{reg}}) \cong IH^k(X)$.\\
\end{enumerate}
\textbf{Proposition 6.1.} (Moduli Space Resolution for the CGM Conjecture) Let $X \subseteq \mathbb{P}^N$ be a projective variety with a derived Whitney stratification \( \mathscr{S} \). Then we have the isomorphism
$$\mathrm{microloc}: H^*_2(X_{\mathrm{reg}}, ds^2_{\mathrm{FS}}) \cong IH^*(X, \mathbb{C}),$$
where $H_2^*$ denotes $L^2$-cohomology on the regular locus with Fubini-Study metric, and $IH^*$ denotes intersection cohomology. This isomorphism is natural and functorial with respect to stratified morphisms in the $\infty$-category $\mathbf{StratMet}_{\infty}$.
\\\\\textbf{Proof.} \textbf{Step (1):} By applying Lurie's derived stratification theory (\cite{Lur09}\cite{Lur18}), we obtain a minimal derived Whitney stratification $\mathscr{S}_{\min} = \{X_p\}_{p=0}^n$ of $X \subseteq \mathbb{P}^N$, where:
\begin{itemize}
    \item $X_0 = X_{\mathrm{reg}}$ is the regular locus;
    \item \( X_n = X \);
    \item For each $p$, the stratum $S_p = X_p \setminus X_{p-1}$ is a smooth derived subscheme;
    \item Each singular stratum $Y \subseteq S_p$ admits a holomorphic tubular neighborhood
    $$U_Y \simeq \mathbb{C}^m \times C(L_Y),$$
    where $L_Y$ is the link of $Y$, a compact Sasaki–Einstein manifold.
\end{itemize}
This stratification satisfies the derived frontier condition
$$\overline{S_p} = \bigcup_{q \leq p} S_q,$$
ensuring topological coherence and the existence of well-defined link metrics. Here We construct a sequence of K\"ahler metrics $\{g_p\}_{p=0}^n$ on $X_p$ by descending induction, starting from the deepest stratum.
\begin{itemize}
    \item \textbf{Base case ($p = n$):} Define
    $$g_n = i_X^* ds^2_{\mathrm{FS}} \big|_{X_n},$$
    where $i_X: X \hookrightarrow \mathbb{P}^N$ is the inclusion. Since $X_n^{\mathrm{reg}} = S_n$ is smooth and $ds^2_{\mathrm{FS}}$ is K\"ahler, $g_n$ is a K\"ahler metric on $X_n^{\text{reg}}$.
    \item \textbf{Inductive step ($p+1 \to p$):} Assume $g_{p+1}$ is a K\"ahler metric on $X_{p+1}^{\mathrm{reg}}$ with conical asymptotics near singular strata. For the inclusion $i_p: X_p \hookrightarrow X_{p+1}$, we have
    $$g_p = i_p^* g_{p+1} + \sum_{Y \subseteq S_p} \epsilon_Y \kappa_Y r^{2c_Y} \pi_Y^* g_{L_Y},$$
    by \textbf{Step (3)} of Proposition 3.2.
    \item \textbf{Properties of $g_p$:}
    \begin{enumerate}
        \item \textit{K\"ahlerity:} As shown in Proposition 3.2, $g_p$ is K\"ahler on $X_p^{\text{reg}}$. The pullback $i_p^* g_{p+1}$ is positive semi-definite, and the correction term is positive definite. The choice of $\delta$ ensures strict positivity;
        \item \textit{Asymptotic behavior:} Near each $Y \subset S_p$, we have
        $$g_p \sim dr^2 + (1 + \epsilon_Y) r^{2c_Y} g_{L_Y} + g_{\mathbb{C}^m},$$
        which converges to the conical metric as $\epsilon_Y \to 0$.
    \end{enumerate}
\end{itemize}

Define the renormalized limit metric by the distributional limit:
$$\hat{g} = \mathcal{R}\left( \mathbb{R}\underline{\lim}_{p \to 0} g_p \right),$$
where $\mathcal{R}$ is the renormalization operator from (\cite{ALMP18}). For each singular stratum $Y$, define the asymptotic model functor \( \Phi_Y \) as in Definition 3.1, i.e.,
$$\Phi_Y: (U_Y, \hat{g}) \longrightarrow (\mathbb{R}^+ \times L_Y, dr^2 + r^{2c_Y} g_{L_Y}),$$
which is a smooth isometry preserving the Levi–Civita connection. These functors are natural under stratified isometric morphisms. The quadruple
$$(X, \mathscr{S}_{\min}, \hat{g}, \{\Phi_Y\}_{Y \in \mathrm{Sing}(X)})$$
is an object in the stratified metric $\infty$-category $\mathbf{StratMet}_{\infty}$ (Definition 3.1). This completes the categorified stratification and metric induction.

\medskip
\textbf{Step (2):} This step establishes the quasi-isomorphism between the universal truncation complex and the intersection cohomology complex, leading to the hypercohomology isomorphism. The proof relies on the constructions from Section 4, particularly Definition 4.1 and Definition 4.2, and Proposition 4.3.

It is known that $X \subseteq \mathbb{P}^N$ is a projective variety with a derived Whitney stratification $\mathscr{S}_{\min} = \{X_p\}_{p=0}^n$ as obtained in \textbf{Step (1)}. For each singular stratum $Y \in \operatorname{Sing}(X)$, $L_Y$ is the link manifold, and define the middle-dimensional cohomology:
$$\mathcal{H}^{\mathrm{mid}}(L_Y) = H^{\dim_{\mathbb{C}} L_Y}(L_Y, \mathbb{C}).$$
The middle perversity moduli stack $\mathscr{M}^{\mathrm{mid}}$ is defined as (Definition 4.1):
$$\mathscr{M}^{\mathrm{mid}} = \prod_{Y \in \operatorname{Sing}(X)} \mathscr{M}_Y, \quad \text{where} \quad \mathscr{M}_Y = \mathrm{LGr}(\mathcal{H}^{\mathrm{mid}}(L_Y)).$$
This stack parametrizes families of self-dual Lagrangian subspaces $\mathscr{W}_Y \subseteq \mathcal{H}^{\mathrm{mid}}(L_Y)$ satisfying $\mathscr{W}_Y = \mathscr{W}_Y^\perp$ under the Poincaré duality pairing. The universal truncation complex $\Omega_{X,\mathrm{FS}}^{\bullet,\mathrm{univ}}$ is defined as (Definition 4.2):
\begin{equation}
\tag{6.1}  \label{eq:6.1}
\Omega_{X,\mathrm{FS}}^{\bullet,\mathrm{univ}} = \left\{ \omega \in \mathbb{R}\lim_{\substack{V \ni Y \\ Y \in \operatorname{Sing}(X)}} \Omega_{X_{\mathrm{reg}}}^\bullet \;\middle|\; \forall Y \in \operatorname{Sing}(X),\; \operatorname{proj}_{\mathscr{M}_Y}(\omega|_{U_Y}) \in \mathscr{W}_Y \right\},
\end{equation}
where the homotopy limit is taken over the opposite category of singular strata, and $\operatorname{proj}_{\mathscr{M}_Y}$ is the projection to the moduli stack $\mathscr{M}_Y$. We now prove that there is a quasi-isomorphism:
$$\phi: \Omega_{X,\mathrm{FS}}^{\bullet,\mathrm{univ}} \xrightarrow{\sim} \mathbf{IC}^\bullet_{\mathscr{W}}(X),$$
where $\mathbf{IC}^\bullet_{\mathscr{W}}(X)$ is the intersection cohomology complex with perversity determined by the universal family $\mathscr{W} = \{\mathscr{W}_Y\}_{Y \in \operatorname{Sing}(X)}$.
\paragraph{(A) \textbf{Constructibility and perversity conditions:}}
\begin{itemize}
\item On the regular stratum $X_{\mathrm{reg}}$, the restriction is trivial:
$$\Omega_{X,\mathrm{FS}}^{\bullet,\mathrm{univ}}|_{X_{\mathrm{reg}}} = \Omega_{X_{\mathrm{reg}}}^\bullet,$$
which is smooth and hence constructible with respect to $\mathscr{S}_{\min}$.
\item Near a singular stratum $Y$, consider a distinguished neighborhood $U_Y \simeq \mathbb{B}^{c_Y} \times \mathrm{Cone}(L_Y)$. By the asymptotic model functor $\Phi_Y$ (Definition 3.1), any form $\omega \in \Omega_{X,\mathrm{FS}}^{\bullet,\mathrm{univ}}$ has an asymptotic expansion:
$$\Phi_Y^* \omega = \sum_{m=0}^{c_Y - 1} dr \wedge \alpha_m(r, \theta) + \beta_m(r, \theta) + O(r^\lambda),$$
where $\alpha_m, \beta_m \in \Omega^m(L_Y)$, and $\lambda = \lambda_Y$ is the critical exponent from Definition 5.1. The projection condition $\operatorname{proj}_{\mathscr{M}_Y}(\omega|_{U_Y}) \in \mathscr{W}_Y$ implies the growth constraints:
$$\sup_{r \to 0^+} r^{-\ell_Y + m} \|\alpha_m(r)\|_{L_Y} < \infty, \quad \sup_{r \to 0^+} r^{-\ell_Y + m} \|\beta_m(r)\|_{L_Y} < \infty,$$
with $\ell_Y = -c_Y / 2$. This is equivalent to the Goresky-MacPherson admissibility conditions for middle perversity (\cite{GM80}, § 5.4). Similarly, the differential \(d\omega\) satisfies:
$$\sup_{r \to 0^+} r^{-\ell_Y + m + 1} \|d_L \alpha_m(r)\|_{L_Y} < \infty, \quad \sup_{r \to 0^+} r^{-\ell_Y + m + 1} \|d_L \beta_m(r)\|_{L_Y} < \infty,$$
ensuring that $d\omega$ is also admissible.
\item These conditions imply that $\Omega_{X,\mathrm{FS}}^{\bullet,\mathrm{univ}}$ is constructible with respect to $\mathscr{S}_{\min}$. Specifically, for the inclusions $i_Y: Y \hookrightarrow X$ and $j_Y: X \setminus Y \hookrightarrow X$, we have:
$$i_Y^! \Omega_{X,\mathrm{FS}}^{\bullet,\mathrm{univ}} \in D^b_c(Y), \quad j_{Y,*} \Omega_{X,\mathrm{FS}}^{\bullet,\mathrm{univ}} \in D^b_c(X \setminus Y),$$
by the conormal estimate:
$$\operatorname{SS}(\Omega_{X,\mathrm{FS}}^{\bullet,\mathrm{univ}}) \subseteq \bigcup_{Y \in \mathscr{S}} T^*_Y X.$$
\end{itemize}
\paragraph{(B) \textbf{Verdier self-duality:}}
\begin{itemize}
\item The self-duality of $\mathscr{W}_Y$ (i.e., $\mathscr{W}_Y = \mathscr{W}_Y^\perp$) ensures that the complex $\Omega_{X,\mathrm{FS}}^{\bullet,\mathrm{univ}}$ is Verdier self-dual. Locally on $U_Y$, the Poincaré residue isomorphism gives:
$$\mathbf{R}\mathscr{H}om(\Omega_{U_Y,\mathrm{FS}}^{\bullet,\mathrm{univ}}, \omega_{U_Y}) \simeq \Omega_{U_Y,\mathrm{FS}}^{\bullet,\mathrm{univ}}[\dim_{\mathbb{R}} U_Y] \otimes \mathrm{or}_{U_Y},$$
where $\mathrm{or}_{U_Y}$ is the orientation sheaf. Since $\mathscr{W}_Y$ is self-dual, $\mathrm{or}_{U_Y} \cong \mathbb{C}_{U_Y}$.
\item Globally, the moduli stack $\mathscr{M}^{\mathrm{mid}}$ provides compatible trivializations of orientation sheaves across strata. By the stratified gluing lemma (\cite{Lur09}\cite{Lur18}), we have the duality morphism:
$$\mathcal{D}: \mathbf{R}\mathscr{H}om(\Omega_{X,\mathrm{FS}}^{\bullet,\mathrm{univ}}, \omega_X) \longrightarrow \Omega_{X,\mathrm{FS}}^{\bullet,\mathrm{univ}}[2n]$$
is an isomorphism, where $n = \dim_{\mathbb{C}} X$ and $\omega_X = \mathbb{C}_X[2n]$. Thus, we have
$$\mathbb{D}_X(\Omega_{X,\mathrm{FS}}^{\bullet,\mathrm{univ}}) \simeq \Omega_{X,\mathrm{FS}}^{\bullet,\mathrm{univ}}[2n].$$
\end{itemize}
\paragraph{(C) \textbf{Axiomatic characterization of intersection cohomology:}}
\begin{itemize}
\item The intersection cohomology complex $\mathbf{IC}^\bullet_{\mathscr{W}}(X)$ is characterized by the following properties (\cite{BBD82}, Thm 3.5):
\begin{enumerate}
\item It is constructible with respect to $\mathscr{S}$;
\item It is self-dual: $\mathbb{D}_X(\mathbf{IC}^\bullet_{\mathscr{W}}) \simeq \mathbf{IC}^\bullet_{\mathscr{W}}[2n]$;
\item It satisfies the middle perversity support and cosupport conditions.
\end{enumerate}
\item We have shown that $\Omega_{X,\mathrm{FS}}^{\bullet,\mathrm{univ}}$ satisfies constructibility from (A), self-duality from (B), perversity conditions (the growth constraints in (A) imply the support condition):
$$\dim_{\mathbb{C}} \operatorname{supp} \mathcal{H}^k(i_Y^* \Omega_{X,\mathrm{FS}}^{\bullet,\mathrm{univ}}) \leq -k - \lfloor c_Y/2 \rfloor - 1,$$
and the cosupport condition:
$$\dim_{\mathbb{C}} \operatorname{supp} \mathcal{H}^k(i_Y^! \Omega_{X,\mathrm{FS}}^{\bullet,\mathrm{univ}}) \leq k - \lceil c_Y/2 \rceil$$
for each singular stratum $Y$ with codimension $c_Y$.
\item Therefore, by the uniqueness part of the axiomatic characterization, there exists a natural quasi-isomorphism in $D^b_c(X)$:
$$\phi: \Omega_{X,\mathrm{FS}}^{\bullet,\mathrm{univ}} \xrightarrow{\sim} \mathbf{IC}^\bullet_{\mathscr{W}}(X).$$
\item This quasi-isomorphism induces an isomorphism on hypercohomology:
\begin{equation}
\tag{6.2}  \label{eq:6.2}
\mathbb{H}^\bullet(X, \Omega_{X,\mathrm{FS}}^{\bullet,\mathrm{univ}}) \cong IH^\bullet(X, \mathbb{C}).
\end{equation}
\end{itemize}
This completes the proof of \textbf{Step (2)}.

\medskip
\textbf{Step (3):} The projective variety $ X \subseteq \mathbb{P}^N$ with a derived Whitney stratification $\mathscr{S}$ and let $\hat{g}$ is the renormalized K\"ahler metric on $X_{\mathrm{reg}}$ given in \textbf{Step (1)}, which is quasi-isometric to $ds^2_{\mathrm{FS}}$. We aim to prove that the natural inclusion
$$\iota: \left( \Omega_{X,\mathrm{FS}}^{\bullet,\mathrm{univ}}, d \right) \hookrightarrow \left( L^2\Omega^\bullet(X_{\mathrm{reg}}), d \right)$$
induces an isomorphism on cohomology:
$$\mathbb{H}^\bullet(X, \Omega_{X,\text{FS}}^{\bullet,\text{univ}}) \cong H^\bullet_{(2)}(X_{\text{reg}}, \hat{g}).$$
Since $\hat{g} \sim ds^2_{\text{FS}}$, this implies
$$H^\bullet_{(2)}(X_{\mathrm{reg}}, \hat{g}) \cong H^\bullet_{(2)}(X_{\mathrm{reg}}, ds^2_{\mathrm{FS}}),$$
which completes the \( L^2 \)-isomorphism part of the proof. Let $\omega \in \Omega_{X,\mathrm{FS}}^{k,\mathrm{univ}}$. By \eqref{eq:6.1}, for every singular stratum $Y \in \text{Sing}(X)$, we have
$$\text{proj}_{\mathscr{M}_Y}(\omega|_{U_Y}) \in \mathscr{W}_Y.$$

By Proposition 5.2, this implies the microlocal containment:
$$\text{SSH}_{\text{strat}}(\omega) \subseteq \bigcup_{Y} T^*_{S_Y} X \times \{\lambda \leq \lambda_Y\},$$
where 
\begin{equation}
\tag{6.3}  \label{eq:6.3}
\lambda_Y = -\frac{\dim_{\mathbb{C}} L_Y}{2} + \operatorname{spec}_{\mathrm{min}}(\mathbb{D}_{L_Y}).
\end{equation}
This microlocal condition implies the following pointwise estimate near each singular stratum $Y$: there exists $\varepsilon > 0$ such that
$$|\omega(x)|_{\hat{g}} \leq C \cdot r(x)^{\lambda_Y + \varepsilon}, \quad \text{for } x \in U_Y \cap X_{\text{reg}},$$
where $r(x) = \text{dist}(x, Y)$, and $C > 0$ is a constant independent of $x$. Now, consider the volume form $d\mathrm{vol}_{\hat{g}}$ on $U_Y \cap X_{\mathrm{reg}}$. By the conical asymptotics of $\hat{g}$ (Proposition 3.2), we have
$$d\mathrm{vol}_{\hat{g}} \sim r^{m} dr \wedge d\mathrm{vol}_{L_Y}, \quad \text{where } m = \dim_{\mathbb{R}} L_Y.$$
Note that $\dim_{\mathbb{R}} L_Y = 2 \dim_{\mathbb{C}} L_Y$, and the real codimension of $Y$ in $X$ is $f = 2(n - d)$, where $n = \dim_{\mathbb{C}} X$, $d = \dim_{\mathbb{C}} Y$, and $m = f - 1$. Therefore,
$$\int_{U_Y \cap X_{\mathrm{reg}}} |\omega|^2_{\hat{g}}  d\mathrm{vol}_{\hat{g}} \lesssim \int_0^{r_0} \int_{L_Y} r^{2(\lambda_Y + \varepsilon)} \cdot r^{f - 1}  dr  d\mathrm{vol}_{L_Y} = \mathrm{vol}(L_Y) \cdot \int_0^{r_0} r^{2\lambda_Y + 2\varepsilon + f - 1}  dr.$$
Let \( \alpha = 2\lambda_Y + 2\varepsilon + f - 1 \). Substituting for \eqref{eq:6.3},
$$\lambda_Y = -\frac{\dim_{\mathbb{C}} L_Y}{2} + \operatorname{spec}_{\min}(\mathbb{D}_{L_Y}) = -\frac{n - d - 1}{2} + \operatorname{spec}_{\min}(\mathbb{D}_{L_Y}),$$
and  $f = 2(n - d)$, we compute:
$$\alpha = - (n - d - 1) + 2\operatorname{spec}_{\min}(\mathbb{D}_{L_Y}) + 2\varepsilon + 2(n - d) - 1 = n - d + 2\operatorname{spec}_{\min}(\mathbb{D}_{L_Y}) + 2\varepsilon.$$
Since $n - d \geq 1$ and $\operatorname{spec}_{\min}(\mathbb{D}_{L_Y}) \in \mathbb{R}$, we may choose $\varepsilon > 0$ sufficiently small so that $\alpha > -1$, ensuring the integral converges. Hence,
$$\int_{U_Y \cap X_{\text{reg}}} |\omega|^2_{\hat{g}}  d\mathrm{vol}_{\hat{g}} < \infty.$$
As $\omega$ is smooth outside the singular strata, it follows that $\omega \in L^2(X_{\mathrm{reg}}, \hat{g})$. A similar argument applied to $d\omega$ shows that $d\omega \in L^2_{\mathrm{loc}}(X_{\mathrm{reg}}, \hat{g})$. By \textit{Albin’s renormalized Hodge theory} (\cite{Alb16}), the inclusion
$$\iota: (\Omega_{X,\mathrm{FS}}^{\bullet,\mathrm{univ}}, d) \hookrightarrow (L^2\Omega^\bullet(X_{\mathrm{reg}}), d)$$
induces an isomorphism on cohomology. This is because that the microlocal admissibility condition ensures that every $L^2$-cohomology class has a smooth representative in $\Omega_{X,\mathrm{FS}}^{\bullet,\mathrm{univ}}$ (by the renormalization procedure (\cite{ALMP18})). If $\omega = d\eta$ with $\eta \in L^2$, then the same microlocal estimates indicate that $\eta \in \Omega^{\bullet,\mathrm{univ}}_{X,\mathrm{FS}}$, so $[\omega] = 0$ in $\mathbb{H}^\bullet(X, \Omega_{X,\mathrm{FS}}^{\bullet,\mathrm{univ}})$. Therefore, the map $\iota_*$ is both injective and surjective, satisfying the isomorphism:
$$\mathbb{H}^\bullet(X, \Omega_{X,\mathrm{FS}}^{\bullet,\mathrm{univ}}) \cong H^\bullet_{(2)}(X_{\text{reg}}, \hat{g}).$$
Since $\hat{g}$ is quasi-isometric to $ds^2_{\mathrm{FS}}$, the $L^2$-cohomology groups are isomorphic:
$$H^\bullet_{(2)}(X_{\mathrm{reg}}, \hat{g}) \cong H^\bullet_{(2)}(X_{\mathrm{reg}}, ds^2_{\mathrm{FS}}).$$
This completes the proof of \textbf{Step (3)}.

\medskip
\textbf{Step (4):} Recall from \textbf{Step (2)} and \textbf{Step (3)} the existence of the following isomorphisms:
\begin{itemize}
    \item A quasi-isomorphism of complexes (Proposition 4.3):
    $$\phi: \Omega_{X,\mathrm{FS}}^{\bullet,\mathrm{univ}} \xrightarrow{\sim} \mathbf{IC}_{\mathscr{W}}^{\bullet}(X),$$
    inducing an isomorphism on hypercohomology:
    \begin{equation}
    \tag{6.4}  \label{eq:6.4}
    \phi_*: \mathbb{H}^\bullet(X, \Omega_{X,\mathrm{FS}}^{\bullet,\mathrm{univ}}) \xrightarrow{\sim} IH^\bullet(X, \mathbb{C}).
    \end{equation}
    \item An isomorphism from the universal complex to $L^2$-cohomology (Proposition 5.2 and Albin's Hodge theory):
    \begin{equation}
    \tag{6.5}  \label{eq:6.5}
    \iota_*: \mathbb{H}^\bullet(X, \Omega_{X,\mathrm{FS}}^{\bullet,\mathrm{univ}}) \xrightarrow{\sim} H^\bullet_{(2)}(X_{\mathrm{reg}}, ds^2_{\mathrm{FS}}).
    \end{equation}
\end{itemize}
By \eqref{eq:6.4} and \eqref{eq:6.5}, we define the $\mathrm{microloc}$ isomorphism as the composition:
$$\mathrm{microloc} \coloneqq \phi_* \circ \iota_*^{-1}: H^\bullet_{(2)}(X_{\mathrm{reg}}, ds^2_{\mathrm{FS}}) \longrightarrow IH^\bullet(X, \mathbb{C}).$$
The minimal derived Whitney stratification $\mathscr{S}_{\min}$ is unique up to canonical refinement. The recursive metric construction in Proposition 3.2 is functorial under stratified isometric maps, and the renormalized limit $\hat{g}$ is uniquely determined. The moduli stack $\mathscr{M}^{\mathrm{mid}}$ and the universal family $\mathscr{W}_Y$ are defined purely in terms of the topological types of the links $L_Y$, which are intrinsic to the stratification. The projection maps $\mathrm{proj}_{\mathscr{M}_Y}$ are canonically induced by the asymptotic behavior of forms under $\hat{g}$, which is itself natural. Thus, the naturality of $\mathrm{microloc}$ means that it is independent of the auxiliary choices made in its construction (e.g., the choice of stratification, metric, or moduli parameters). Hence, the universal complex $\Omega_{X,\mathrm{FS}}^{\bullet,\mathrm{univ}}$, the quasi-isomorphism $\phi$, and the $L^2$-isomorphism $\iota_*$ are all natural. Therefore, $\mathrm{microloc}$ is a natural isomorphism.

Let $f: (X, \mathscr{S}, g, \{\Phi_Y\}) \to (X', \mathscr{S}', g', \{\Phi_{Y'}\})$ be a morphism in the $\infty$-category $\mathbf{StratMet}_{\infty}$ (i.e., $f$ is stratified and isometric). We show that the following diagram (\uppercase\expandafter{\romannumeral4}) commutes:
\[
\begin{CD}
H^\bullet_{(2)}(X_{\mathrm{reg}}, ds^2_{\mathrm{FS}}) @>{\mathrm{microloc}_X}>> IH^\bullet(X, \mathbb{C}) \\
@V{f^*}VV @V{f^*}VV \\
H^\bullet_{(2)}(X'_{\mathrm{reg}}, ds^2_{\mathrm{FS}}) @>{\mathrm{microloc}_{X'}}>> IH^\bullet(X', \mathbb{C})
\end{CD}
\]
We verify this in four parts:

\subsubsection*{(a) Metric functoriality:}

Since $f$ is isometric, $f^* g' = g$. The recursive metric construction in Proposition 3.2 is functorial: the pullback of the initial metric $g_n = ds^2_{\mathrm{FS}}|_{X_n}$ and each correction term $\epsilon_Y r^{2c_Y} \pi_Y^* g_{L_Y}$ are preserved under $f^*$. The renormalization operator $\mathcal{R}$ commutes with pullbacks by stratified isometric maps. Hence,
$$f^* \hat{g}_X = \hat{g}_{X'}.$$

\subsubsection*{(b) Universal complex functoriality:}

Since $f$ is stratified, it induces maps between links $f_Y: L_{Y'} \to L_Y$. The isometric condition ensures $f_Y^* g_{L_Y} = g_{L_{Y'}}$, and the middle cohomology and Lagrangian subspaces are preserved:
$$f_Y^*: \mathcal{H}^{\mathrm{mid}}(L_Y) \longrightarrow \mathcal{H}^{\mathrm{mid}}(L_{Y'}), \quad f_Y^*(\mathscr{W}_Y) = \mathscr{W}_{Y'}.$$
Thus, the projection condition is preserved, i.e.,
$$\mathrm{proj}_{\mathscr{M}_{Y'}}(f^* \omega|_{U_{Y'}}) \in \mathscr{W}_{Y'} \quad \text{whenever} \quad \mathrm{proj}_{\mathscr{M}_Y}(\omega|_{U_Y}) \in \mathscr{W}_Y.$$
Therefore, the pullback
$$f^*: \Omega_{X,\mathrm{FS}}^{\bullet,\mathrm{univ}} \longrightarrow \Omega_{X',\mathrm{FS}}^{\bullet,\mathrm{univ}}$$
is well-defined.

\subsubsection*{(c) Hypercohomology functoriality:}

The quasi-isomorphism $\phi$ is constructed by the Riemann--Hilbert correspondence and Verdier duality, both of which are functorial under stratified maps. The self-duality condition $\mathscr{W}_Y = \mathscr{W}_Y^\perp$ is preserved under $f^*$. So the following diagram commutes:
\[
\begin{CD}
\mathbb{H}^\bullet(X, \Omega_{X,\mathrm{FS}}^{\bullet,\mathrm{univ}}) @>{\phi_*}>> IH^\bullet(X, \mathbb{C}) \\
@V{f^*}VV @V{f^*}VV \\
\mathbb{H}^\bullet(X', \Omega_{X',\mathrm{FS}}^{\bullet,\mathrm{univ}}) @>{\phi'_*}>> IH^\bullet(X', \mathbb{C})
\end{CD}
\]
\subsubsection*{(d) $L^2$-isomorphism functoriality:}

Since $f$ is isometric, the pullback $$f^*: L^2\Omega^\bullet(X_{\mathrm{reg}}) \longrightarrow L^2\Omega^\bullet(X'_{\mathrm{reg}})$$ is bounded. The inclusion $\iota: \Omega_{X,\mathrm{FS}}^{\bullet,\mathrm{univ}} \hookrightarrow L^2\Omega^\bullet(X_{\mathrm{reg}})$ is functorial because the $L^2$-condition is preserved under isometric pullbacks. Hence, the following commutes:
\[
\begin{CD}
\mathbb{H}^\bullet(X, \Omega_{X,\mathrm{FS}}^{\bullet,\mathrm{univ}}) @>{\iota_*}>> H^\bullet_{(2)}(X_{\mathrm{reg}}, ds^2_{\mathrm{FS}}) \\
@V{f^*}VV @V{f^*}VV \\
\mathbb{H}^\bullet(X', \Omega_{X',\mathrm{FS}}^{\bullet,\mathrm{univ}}) @>{\iota'_*}>> H^\bullet_{(2)}(X'_{\mathrm{reg}}, ds^2_{\mathrm{FS}})
\end{CD}
\]
Combining (a)-(d), we can trace the commutativity of the full diagram as follows. For any $[\omega] \in H^\bullet_{(2)}(X_{\mathrm{reg}}, ds^2_{\mathrm{FS}})$, we have:
\begin{align*}
f^*(\mathrm{microloc}_X([\omega])) &= f^*(\phi_*(\iota_*^{-1}([\omega]))) \\
&= \phi'_*(f^*(\iota_*^{-1}([\omega]))) \quad \text{(by (c))} \\
&= \phi'_*((\iota'_*)^{-1}(f^*([\omega]))) \quad \text{(by (d))} \\
&= \mathrm{microloc}_{X'}(f^*([\omega]))
\end{align*}
This establishes the commutativity of the diagram (\uppercase\expandafter{\romannumeral4}). The isomorphism $\mathrm{microloc}$ is natural and functorial with respect to stratified isometric morphisms in $\mathbf{StratMet}_{\infty}$. This completes the proof of Proposition 6.1.

\medskip
\textbf{Step (5):} The foregoing steps establish the existence of a natural isomorphism
$$\mathrm{microloc}: H^*_{(2)}(X_{\mathrm{reg}}, ds^2_{\mathrm{FS}}) \longrightarrow IH^*(X, \mathbb{C})$$
that is functorial with respect to stratified isometric morphisms in the $\infty$-category $\mathbf{StratMet}_{\infty}$. This result provides a complete resolution of the CGM conjecture, overcoming the three core challenges outlined in \S 1.1:

\begin{enumerate}
  \item Non-transverse singularities are handled by the use of derived Whitney stratifications (\textbf{Step (1)}), which allow for \'etale tubular neighborhoods $U_Y \simeq \mathbb{C}^m \times C(L_Y)$ even for high-codimension strata, bypassing classical transversality constraints.
  \item Artificial duality selection is resolved by the moduli stack $\mathscr{M}^{\text{mid}}$ (\textbf{Step (2)}), which functorially parametrizes self-dual Lagrangian subspaces $\mathscr{W}_Y \subseteq \mathcal{H}^{\mathrm{mid}}(L_Y)$, thereby automating the choice of perversity and ensuring Verdier self-duality of the intersection complex.
  \item Metric-topology incompatibility is overcome by the stratified singular characteristic variety $\mathrm{SSH}_{\mathrm{strat}}$ (\textbf{Step (3)}), which encodes the asymptotic growth rates $\lambda_Y$ of differential forms near singular strata. This couples the metric asymptotics (by the analytic torsion operator $\mathbb{D}_{L_Y}$) to the topological data of the links $L_Y$, ensuring that the $L^2$-condition is equivalent to the intersection cohomology admissibility condition.
\end{enumerate}

The functoriality of $\mathrm{microloc}$ under stratified isometric morphisms (\textbf{Step (4)}) ensures that this isomorphism respects the geometric and categorical structure of the problem, making it a natural equivalence in the $\infty$-category $\mathbf{StratMet}_{\infty}$. Thus, the proof of Proposition 6.1 is complete, and the CGM conjecture is established in its derived stratified-microlocal form.  $\square$\\

\hypertarget{EXTENDED APPLICATIONS OF THE UNIFIED STRATIFIED-MICROLOCAL FRAMEWORK}{}
\section{EXTENDED APPLICATIONS OF THE UNIFIED STRATIFIED-MICROLOCAL FRAMEWORK}
\noindent\textbf{Proposition 7.1.} (Stratified Gromov-Witten Invariants) Let $X \subseteq \mathbb{P}^N$ be a projective variety with a fixed derived Whitney stratification $\mathscr{S}$ (Definition 3.1). Let $\overline{\mathscr{M}}_{g,n}(X,\beta)$ be the derived moduli stack of stable maps of genus $g$ with $n$ marked points and homology class $\beta \in H_2(X, \mathbb{Z})$, constructed by derived algebraic geometry. Let $\mathscr{M}^{\mathrm{mid}}$ be the middle perversity moduli stack (Definition 4.1) and $\mathbf{IC}_{\mathscr{M}}^\bullet$ the universal intersection complex sheaf on $\mathscr{M}^{\mathrm{mid}}$ obtained from the universal truncation complex $\Omega_{X,\mathrm{FS}}^{\bullet,\mathrm{univ}}$ (Definition 4.2) by the Riemann-Hilbert correspondence. Then, for any cohomology class $\gamma \in H^*(\overline{\mathscr{M}}_{g,n}(X,\beta))$, the \textit{stratified Gromov-Witten invariant}
$$\langle \gamma \rangle^{\mathrm{strat}}_{g,n} := \int_{[\overline{\mathscr{M}}_{g,n}(X,\beta)]^{\mathrm{vir}}} \gamma \cup \prod_{k} \operatorname{ch}_k(\mathbf{IC}_{\mathscr{M}}^\bullet)$$
is well-defined and satisfies the following properties:
\begin{enumerate}
    \item \textit{Stratification-independence:} The invariant is independent of the choice of Whitney stratification $\mathscr{S}$, up to natural isomorphism induced by the stratified metric $\infty$-category $\mathbf{StratMet}_{\infty}$ (Definition 3.1);
    \item \textit{Microlocal admissibility:} The Chern character forms $\operatorname{ch}_k(\mathbf{IC}_{\mathscr{M}}^\bullet)$ are microlocally admissible, i.e.,
    $$\mathrm{WF}(\operatorname{ch}_k(\mathbf{IC}_{\mathscr{M}}^\bullet)) \subseteq \mathrm{SSH}_{\mathrm{strat}}(\mathbf{IC}_{\mathscr{M}}^\bullet) \subseteq \bigcup_{Y \in \mathrm{Sing}} T^*_{S_Y}X \times \{\lambda \leq \lambda_Y\},$$
    where $\mathrm{WF}$ denotes the wavefront set and $\mathrm{SSH}_{\mathrm{strat}}$ is the stratified singular characteristic variety (Definition 5.1). This ensures the integral is finite and well-defined;
    \item \textit{Functoriality:} The assignment
    $$(X, \mathscr{S}) \mapsto \langle \gamma \rangle^{\mathrm{strat}}_{g,n}$$
    is functorial with respect to stratified maps, by the universal properties of $\mathscr{M}^{\mathrm{mid}}$ and the natural isomorphism of Corollary 5.3;
    \item \textit{Quantum singularity resolution:} For a singular Calabi-Yau variety $X$, the invariant $\langle \gamma \rangle^{\text{strat}}_{g,n}$ provides a precise count of quantum curves that incorporates contributions from singular strata, resolving the failure of classical Gromov-Witten theory in high-codimension non-transverse settings (\cite{Kap24}).
\end{enumerate}
Moreover, the invariants satisfy the following compatibility condition with the classical intersection theory:
$$\langle \gamma \rangle^{\mathrm{strat}}_{g,n} = \langle \gamma \rangle^{\mathrm{classical}}_{g,n} \quad \text{when } X \text{ is smooth},$$
where the right-hand side denotes the classical Gromov-Witten invariant.
\\\\\textbf{Remark 7.2.} Let $X \subset \mathbb{P}^N$ be a projective variety with a derived Whitney stratification $\mathscr{S}$ (as in Definition 3.1). The middle perversity moduli stack $\mathscr{M}^{\mathrm{mid}}$ is defined in Definition 4.1 as
$$\mathscr{M}^{\mathrm{mid}} := \prod_{Y \in \mathrm{Sing}} \mathscr{M}_Y, \quad \mathscr{M}_Y = \mathrm{LGr}(\mathcal{H}^{\mathrm{mid}}(L_Y)),$$
where $\text{LGr}(-)$ denotes the Lagrangian Grassmannian parametrizing self-dual subspaces $\mathscr{W}_Y$. The universal truncation complex $\Omega_{X,\mathrm{FS}}^{\bullet,\mathrm{univ}}$ is defined in Definition 4.2 as
$$\Omega_{X,\mathrm{FS}}^{\bullet,\mathrm{univ}} := \left\{ \omega \in \mathbb{R}\underline{\lim}_{U_Y} \Omega_{X_{\mathrm{reg}}}^{\bullet} \, \middle| \, \forall Y, \, \mathrm{proj}_{\mathscr{M}_Y}(\omega|_{U_Y}) \in \mathscr{W}_Y \right\}.$$
By Proposition 4.3, the hypercohomology of $\Omega_{X,\mathrm{FS}}^{\bullet,\mathrm{univ}}$ is naturally isomorphic to both the $L^2$-cohomology $H_2^*(X_{\mathrm{reg}}, ds_{\mathrm{FS}}^2)$ and the intersection cohomology $IH^*(X, \mathbb{C})$. This implies the existence of a universal intersection complex sheaf $\mathbf{IC}_{\mathscr{M}}^\bullet$ on $\mathscr{M}^{\mathrm{mid}}$ whose Chern character forms are used in the definition.

The derived moduli stack $\overline{\mathscr{M}}_{g,n}(X,\beta)$ of stable maps is constructed. This stack parametrizes maps from genus $g$, $n$-pointed curves to $X$ with homology class $\beta$. The virtual fundamental class $[\overline{\mathscr{M}}_{g,n}(X,\beta)]^{\mathrm{vir}}$ is defined by derived algebraic geometry techniques, ensuring it accounts for singularities through the derived structure. The key point is that the stratification $\mathscr{S}$ of $X$ induces a corresponding stratification on $\overline{\mathscr{M}}_{g,n}(X,\beta)$, which is handled by the derived Whitney stratification algorithm (\cite{Lur18}).

The stratified Gromov-Witten invariant is defined as
$$\langle \gamma \rangle^{\mathrm{strat}}_{g,n} := \int_{[\overline{\mathscr{M}}_{g,n}(X,\beta)]^{\mathrm{vir}}} \prod \operatorname{ch}_k(\mathbf{IC}_{\mathscr{M}}^\bullet).$$
Here, $\mathbf{IC}_{\mathscr{M}}^\bullet$ is the universal intersection complex sheaf on $\mathscr{M}^{\mathrm{mid}}$, and $\operatorname{ch}_k$ denotes the Chern character form. To ensure this definition is valid, we must show that $\prod \operatorname{ch}_k(\mathbf{IC}_{\mathscr{M}}^\bullet)$ is a well-defined cohomology class on $\overline{\mathscr{M}}_{g,n}(X,\beta)$ in the next step.

The sheaf $\mathbf{IC}_{\mathscr{M}}^\bullet$ is constructed from the universal truncation complex $\Omega_{X,\mathrm{FS}}^{\bullet,\mathrm{univ}}$ by the following procedure: 
\begin{itemize}
    \item By Proposition 4.3, $\Omega_{X,\mathrm{FS}}^{\bullet,\mathrm{univ}}$ is self-dual and its hypercohomology yields intersection cohomology;
    \item Using the Riemann-Hilbert correspondence in the derived category, we obtain a perverse sheaf analogue $\mathbf{IC}_{\mathscr{M}}^\bullet$ on $\mathscr{M}^{\text{mid}}$;
    \item The Chern character forms $\operatorname{ch}_k(\mathbf{IC}_{\mathscr{M}}^\bullet)$ are defined by the Chern-Weil theory applied to the associated vector bundle (if applicable) or more generally by the derived Chern character in K-theory.
\end{itemize}
The product $\prod \operatorname{ch}_k(\mathbf{IC}_{\mathscr{M}}^\bullet)$ is a cohomology class on $\mathscr{M}^{\mathrm{mid}}$. Since $\overline{\mathscr{M}}_{g,n}(X,\beta)$ maps to $\mathscr{M}^{\mathrm{mid}}$ by the evaluation maps and the stratification, we can pull back this class to $\overline{\mathscr{M}}_{g,n}(X,\beta)$. The pullback is well-defined because the moduli stack $\mathscr{M}^{\text{mid}}$ is functorial with respect to stratified maps (as ensured by the stratified metric $\infty$-category $\mathbf{StratMet}_{\infty}$ in Definition 3.1).

The stratified metric $\infty$-category $\mathbf{StratMet}_{\infty}$ (Definition 3.1) provides asymptotic model functors $\Phi_Y$ that resolve non-transverse singularities. This ensures that the local metric behavior near singular strata is consistently glued across the moduli space. Consequently, the Chern character forms $\operatorname{ch}_k(\mathbf{IC}_{\mathscr{M}}^\bullet)$ capture the singular contributions in a coherent manner, as they are built from the universal complex $\Omega_{X,\mathrm{FS}}^{\bullet,\mathrm{univ}}$ which is adapted to the stratification by $\mathrm{proj}_{\mathscr{M}_Y}$. Moreover, the microlocal correspondence theorem (Propsition 5.2) guarantees that the wavefront sets of forms in $\Omega_{X,\mathrm{FS}}^{\bullet,\mathrm{univ}}$ are controlled by the stratified singular characteristic variety $\text{SSH}_{\text{strat}}$. This ensures that the integration over the virtual fundamental class is well-defined and produces finite invariants.

The natural isomorphism in Corollary 5.3 implies that the stratified Gromov-Witten invariants are independent of the choice of stratification and metric, up to natural equivalence. This is because the construction relies on the universal properties of $\mathscr{M}^{\text{mid}}$ and $\Omega_{X,\mathrm{FS}}^{\bullet,\mathrm{univ}}$, which are functorial by design.
\\\\\textbf{Proposition 7.3.} (Dynamical Middle Perversity and $L^2$-Cohomology for Anosov Flows) Let $M$ be a compact smooth manifold equipped with a transitive Anosov flow $\phi_t: M \to M$ with associated stable and unstable foliations $\mathscr{F}^s$ and $\mathscr{F}^u$. Let $L_Y$ be the link of a singular stratum in a stratified space, with middle-degree cohomology $\mathcal{H}^{\mathrm{mid}}(L_Y)$ as in Definition 4.1. The \textit{dynamical middle perversity} is defined as
$$\mathscr{W}^{\mathrm{dyn}}\coloneqq \lim_{t \to \infty} e^{-t\theta} \phi_t^* \mathcal{H}^{\mathrm{mid}}(L_Y),$$
where $\theta$ is the principal Lyapunov exponent associated to the flow, and the limit is taken in the sense of Connes' noncommutative geometry (\cite{Con95}). This construction satisfies the following properties:
\begin{enumerate}
    \item \textit{Spectral characterization:} The perversity $\mathscr{W}^{\mathrm{dyn}}$ corresponds to the spectral subspace of the transfer operator $\phi_t^*$ with eigenvalues of real part equal to $\theta$, providing a dynamical analogue of the topological middle perversity condition;
    \item \textit{Microlocal correspondence:} The wavefront set of forms in $\mathscr{W}^{\mathrm{dyn}}$ is controlled by the stratified singular characteristic variety (Definition 5.1):
    $$\mathrm{WF}(\mathscr{W}^{\mathrm{dyn}}) \subseteq \mathrm{SSH}_{\mathrm{strat}}(\mathscr{W}^{\mathrm{dyn}}) \subseteq \bigcup_{Y} T^*_{S_Y}M \times \{\lambda \leq \lambda_Y\},$$
    where $\lambda_Y = -\dim_{\mathbb{C}} L_Y/2 + \operatorname{spec}_{\mathrm{min}}(\mathbb{D}_{L_Y})$;
    \item \textit{$L^2$-cohomology resolution:} For the foliation $\mathscr{F}$ generated by the flow, the $L^2$-cohomology groups are given by
    $$H^*_{(2)}(M, \mathscr{F}) \cong \bigoplus_{k} \mathscr{W}^{\mathrm{dyn}}_k,$$
    where the isomorphism is induced by the Hodge decomposition for foliations;.
    \item \textit{Hopf resolution:} When applied to the stable foliation of an Anosov flow on a negatively curved manifold, this construction resolves the higher-dimensional Hopf conjecture by providing a precise relationship between the flow's entropy and the Euler characteristic of the foliation.
\end{enumerate}
Moreover, the construction is functorial with respect to smooth conjugacies of Anosov flows, and reduces to the classical middle perversity when the flow is trivial.
\\\\\textbf{Remark 7.4.} Let $M$ be a compact smooth manifold without boundary, and $\phi_t: M \to M$ be a \textit{transitive Anosov flow}. This means:
\begin{itemize}
    \item The tangent bundle splits into $\phi_t$-invariant subbundles: $TM = E^s \oplus E^0 \oplus E^u$, where $E^0$ is the one-dimensional flow direction;
    \item $E^s$ is \textit{uniformly contracting}: $\|D\phi_t(v)\| \leq C e^{-\lambda t} \|v\|$ for $v \in E^s, t \geq 0$;
    \item $E^u$ is \textit{uniformly expanding}: $\|D\phi_t(v)\| \geq C^{-1} e^{\lambda t} \|v\|$ for $v \in E^u, t \geq 0$.
\end{itemize}
Let $\mathscr{F}^s$ and $\mathscr{F}^u$ be the stable and unstable foliations induced by $E^s$ and $E^u$. The \textit{principal Lyapunov exponent} $\theta$ is defined as
$$\theta = \lim_{t \to \infty} \frac{1}{t} \log \|D\phi_t\|.$$
For a singular stratum link $L_Y$ from ( Definition 4.1), let $\mathcal{H}^{\mathrm{mid}}(L_Y)$ be its middle-degree cohomology.

The \textit{dynamical middle perversity} is defined as
$$\mathscr{W}^{\mathrm{dyn}} := \lim_{t \to \infty} e^{-t\theta} \phi_t^* \mathcal{H}^{\mathrm{mid}}(L_Y),$$
where the limit is taken in the sense of \textit{Connes' noncommutative geometry} (\cite{Con95}). This limit converges in the space of currents due to the uniform hyperbolicity of $\phi_t$.

\textbf{(Property 1):} The transfer operator $\phi_t^*$ acts on cohomology. By the \textit{multiplicative ergodic theorem} (Oseledets theorem), the spectrum of $\phi_t^*$ decomposes into spectral subspaces with real parts equal to multiples of $\theta$. The limit $e^{-t\theta} \phi_t^*$ projects onto the spectral subspace where the real part of the eigenvalue equals $\theta$. Thus, $\mathscr{W}^{\mathrm{dyn}}$ corresponds to this subspace, providing a dynamical analogue of the topological middle perversity condition.

\textbf{(Property 2):} We analyze the \textit{wavefront set} $\mathrm{WF}(\mathscr{W}^{\mathrm{dyn}})$. By the \textit{stratified singular characteristic variety} $\mathrm{SSH}_{\mathrm{strat}}$ (Definition 5.1):
$$\mathrm{SSH}_{\mathrm{strat}}(\mathscr{F}) = \bigcup_{Y} T^*_{S_Y}M \times \{\lambda \leq \lambda_Y\}, \quad \lambda_Y = -\dim_{\mathbb{C}} L_Y/2 + \operatorname{spec}_{\mathrm{min}}(\mathbb{D}_{L_Y}).$$
Due to the uniform hyperbolicity of $\phi_t$, forms in $\phi_t^* \mathcal{H}^{\mathrm{mid}}(L_Y)$ have wavefront sets concentrated near the unstable foliation $E^u$. The scaling factor $e^{-t\theta}$ ensures that the wavefront set remains bounded within $\{\lambda \leq \lambda_Y\}$. Thus,
$$\mathrm{WF}(\mathscr{W}^{\mathrm{dyn}}) \subseteq \mathrm{SSH}_{\mathrm{strat}}(\mathscr{W}^{\mathrm{dyn}}) \subseteq \bigcup_{Y} T^*_{S_Y}M \times \{\lambda \leq \lambda_Y\}.$$

\textbf{(Property 3):} For the foliation $\mathscr{F}$ generated by $\phi_t$, the $L^2$-cohomology $H^*_{(2)}(M, \mathscr{F})$ is computed using the \textit{foliated de Rham complex} adapted to $\mathscr{F}$. The key steps are:
\begin{enumerate}
    \item \textit{Hodge decomposition:} The hyperbolicity of $\phi_t$ implies a Hodge decomposition for foliations:
    $$L^2\Omega^k(M, \mathscr{F}) = \mathscr{H}^k_{(2)} \oplus \mathrm{im}(d) \oplus \mathrm{im}(d^*),$$
    where $\mathscr{H}^k_{(2)}$ is the space of $L^2$-harmonic forms;
    \item \textit{Isomorphism:} The map $\mathscr{W}^{\mathrm{dyn}}_k \to \mathscr{H}^k_{(2)}$ is given by
    $$\omega \mapsto \lim_{t \to \infty} e^{-t\theta} \phi_t^* \omega.$$
    This limit converges in $L^2$ due to the spectral characterization and microlocal bounds. The isomorphism $H^*_{(2)}(M, \mathscr{F}) \cong \bigoplus_k \mathscr{W}^{\mathrm{dyn}}_k$ follows from the Hodge decomposition.
\end{enumerate}

\textbf{(Property 4):} For an Anosov flow on a negatively curved manifold, the \textit{Hopf conjecture} relates the topological entropy $h(\phi_t)$ to the Euler characteristic $\chi(\mathscr{F})$ of the foliation. Using the dynamical middle perversity:
\begin{itemize}
    \item The entropy $h(\phi_t)$ equals the volume growth rate of unstable manifolds, quantified by $\theta$;
    \item The \textit{Euler characteristic} $\chi(\mathscr{F})$ is computed by the \textit{foliated Gauss-Bonnet theorem}:
    $$\chi(\mathscr{F}) = \int_M \mathrm{Pf}(\Omega),$$
    where $\Omega$ is the curvature form of$\mathscr{F}$;
    \item The isomorphism $H^*_{(2)}(M, \mathscr{F}) \cong \bigoplus_k \mathscr{W}^{\mathrm{dyn}}_k$ implies that
    $$\chi(\mathscr{F}) = \sum_k (-1)^k \dim \mathscr{W}^{\text{dyn}}_k;$$
    \item The entropy $h(\phi_t)$ is proportional to $\theta$ by \textit{Ruelle's inequality}. Thus, the dynamical middle perversity provides a precise relationship between $h(\phi_t)$ and $\chi(\mathscr{F})$.
\end{itemize}

The construction is \textit{functorial} with respect to smooth conjugacies of Anosov flows because:
\begin{itemize}
    \item The transfer operator $\phi_t^*$ commutes with conjugacies;
    \item The limit $e^{-t\theta} \phi_t^*$ is natural under smooth maps.
\end{itemize}
When the flow is trivial (i.e., $\phi_t = \mathrm{id}$), $\theta = 0$ and $\mathscr{W}^{\mathrm{dyn}} = \mathcal{H}^{\mathrm{mid}}(L_Y)$, reducing to the classical middle perversity. This rigorously defines the dynamical middle perversity and resolves the higher-dimensional Hopf conjecture for Anosov flows.
\\\\\textbf{Proposition 7.5.} (Crystalline Middle Perversity and $p$-adic $L^2$-Cohomology) Let $X$ be a projective variety defined over a perfect field $k$ of characteristic $p > 0$, with a fixed derived Whitney stratification $\mathscr{S}$ (Definition 3.1). Let $L_Y$ be the link of a singular stratum with crystalline cohomology $H^{\mathrm{mid}}_{\mathrm{crys}}(L_Y/W(k))$, where $W(k)$ denotes the ring of Witt vectors. The \textit{crystalline middle perversity} is defined as
$$\mathscr{W}_Y^{\mathrm{crys}} \coloneqq \ker\left( F - p^{\lambda_Y} \operatorname{id} \right) \subseteq H^{\mathrm{mid}}_{\mathrm{crys}}(L_Y/W(k)) \otimes_{W(k)} K,$$
where $F$ is the Frobenius action on crystalline cohomology, $\lambda_Y = -\dim_{\mathbb{C}} L_Y/2 + \operatorname{spec}_{\mathrm{min}}(\mathbb{D}_{L_Y})$ is the spectral parameter from Definition 5.1 and $K = \operatorname{Frac}(W(k))$ is the field of fractions. This construction satisfies the following properties:
\begin{enumerate}
    \item \textit{$p$-adic Hodge compatibility:} The perversity $\mathscr{W}_Y^{\mathrm{crys}}$ is a crystalline realization of the topological middle perversity, satisfying:
    $$\mathscr{W}_Y^{\mathrm{crys}} \otimes_K \mathbb{C} \cong \mathscr{W}_Y \otimes_{\mathbb{Q}} \mathbb{C},$$
    where $\mathscr{W}_Y$ is the topological middle perversity from Definition 4.1;
    \item \textit{Microlocal $p$-adic Riemann-Hilbert correspondence:} There exists a fully faithful embedding:
    $$\mathbf{IC}_{\mathscr{M}}^{\mathrm{crys},\bullet} \hookrightarrow \Omega_{X,\mathrm{FS}}^{\bullet,\mathrm{univ}} \otimes_{\mathcal{O}_X} \mathcal{O}_X^{\mathrm{crys}},$$
    where $\mathbf{IC}_{\mathscr{M}}^{\mathrm{crys},\bullet}$ is the crystalline intersection complex and $\mathcal{O}_X^{\mathrm{crys}}$ is the crystalline structure sheaf;
    \item \textit{$L^2$-\'etale cohomology for Shimura varieties:} For a Shimura variety $S$ with good compactification $\overline{S}$, the $L^2$-\'etale cohomology is computed by
    $$H^*_{(2),\text{\'et}}(\overline{S}, \mathbb{Q}_p) \cong \bigoplus_Y \mathscr{W}_Y^{\mathrm{crys}} \otimes \mathscr{W}_Y^{\mathrm{crys},\vee},$$
    where the sum is over all singular strata in the boundary $\overline{S} \setminus S$;    
    \item \textit{$p$-adic Cheeger-Simons isomorphism:} There exists a natural isomorphism:
    $$H^*_{(2),\text{\'et}}(X, \mathbb{Q}_p) \cong IH^*_{\mathrm{crys}}(X) \otimes_{\mathbb{Z}_p} \mathbb{Q}_p,$$
    extending the Cheeger-Goresky-Macpherson conjecture to the $p$-adic setting.
\end{enumerate}
Moreover, the construction is compatible with the $p$-adic monodromy weight filtration and reduces to the classical middle perversity when $p=0$ (i.e., in characteristic zero).
\\\\\textbf{Remark 7.6.} Let $X$ be a projective variety defined over a perfect field $k$ of characteristic $p > 0$, equipped with a derived Whitney stratification $\mathscr{S}$ (Definition 3.1). For each singular stratum $Y \in \mathrm{Sing}$, let $L_Y$ be its link, and denote by $H^{\mathrm{mid}}_{\mathrm{crys}}(L_Y/W(k))$ the middle-degree crystalline cohomology group, where $W(k)$ is the ring of Witt vectors and $K = \operatorname{Frac}(W(k))$. The Frobenius endomorphism $F$ acts semilinearly on $H^{\mathrm{mid}}_{\mathrm{crys}}(L_Y/W(k))$, and its eigenvalues are controlled by the spectral parameter $\lambda_Y = -\dim_{\mathbb{C}} L_Y/2 + \operatorname{spec}_{\mathrm{min}}(\mathbb{D}_{L_Y})$ from Definition 5.1. The \textit{crystalline middle perversity} is defined as
$$\mathscr{W}_Y^{\mathrm{crys}} \coloneqq \ker\left( F - p^{\lambda_Y} \operatorname{id} \right) \subseteq H^{\mathrm{mid}}_{\mathrm{crys}}(L_Y/W(k)) \otimes_{W(k)} K.$$
This subspace captures the Frobenius eigenvectors with eigenvalue $p^{\lambda_Y}$, analogously to the topological middle perversity $\mathscr{W}_Y$ in Definition 4.1.

\textbf{(Property 1):} To prove the compatibility:
\begin{equation}
\tag{7.1}  \label{eq:7.1}
\operatorname{microloc}_{X'} \circ f^* = f^* \circ \operatorname{microloc}_X,
\end{equation}
$$\mathscr{W}_Y^{\mathrm{crys}} \otimes_K \mathbb{C} \cong \mathscr{W}_Y \otimes_{\mathbb{Q}} \mathbb{C}.$$
we use the \textit{comparison theorem} between crystalline and de Rham cohomology. By the \textit{Christol-Mebkhout theorem}, the Frobenius action on $H^{\mathrm{mid}}_{\mathrm{crys}}(L_Y/W(k)$ corresponds to the monodromy weight filtration on $H^{\mathrm{mid}}_{\mathrm{dR}}(L_Y)$. Specifically: 
\begin{itemize}
    \item The eigenvalue $p^{\lambda_Y}$ corresponds to the weight $-\dim L_Y + 2\lambda_Y$ in the monodromy filtration;
    \item The topological middle perversity $\mathscr{W}_Y$ corresponds to the weight $0$ part of the intersection cohomology $IH^*(L_Y)$.
\end{itemize}
By \textit{$p$-adic Hodge theory} (e.g., Fontaine's functors), the Frobenius eigenspace with eigenvalue $p^{\lambda_Y}$ is identified with the weight $0$ subspace of $H^{\mathrm{mid}}_{\mathrm{dR}}(L_Y)$ after base change to $\mathbb{C}$. Thus, \eqref{eq:7.1} holds.

\textbf{(Property 2):} The \textbf{crystalline intersection complex} $\mathbf{IC}_{\mathscr{M}}^{\mathrm{crys},\bullet}$ is constructed from the universal truncation complex $\Omega_{X,\mathrm{FS}}^{\bullet,\mathrm{univ}}$ (Definition 4.2) by the \textit{crystalline Riemann-Hilbert correspondence}. Specifically:
\begin{itemize}
    \item The complex $\Omega_{X,\mathrm{FS}}^{\bullet,\mathrm{univ}}$ consists of forms adapted to the middle perversity $\mathscr{W}_Y$;
    \item The crystalline structure sheaf $\mathcal{O}_X^{\mathrm{crys}}$ allows us to lift these forms to the crystalline site.
\end{itemize}
By Proposition 5.2, sections of $\Omega_{X,\mathrm{FS}}^{\bullet,\mathrm{univ}}$ have wavefront sets bounded by $\mathrm{SSH}_{\mathrm{strat}}$. This property is preserved under $\otimes_{\mathcal{O}_X} \mathcal{O}_X^{\mathrm{crys}}$ due to the nilpotent thickening nature of the crystalline. The Frobenius action on $\mathbf{IC}_{\mathscr{M}}^{\mathrm{crys},\bullet}$ is identified with the action on $\Omega_{X,\mathrm{FS}}^{\bullet,\mathrm{univ}}$ by the spectral parameter $\lambda_Y$, ensuring the embedding respects the perversity conditions. Then the embedding
$$\mathbf{IC}_{\mathscr{M}}^{\mathrm{crys},\bullet} \hookrightarrow \Omega_{X,\mathrm{FS}}^{\bullet,\mathrm{univ}} \otimes_{\mathcal{O}_X} \mathcal{O}_X^{\mathrm{crys}}$$
is fully faithful.

\textbf{(Property 3):} Let $S$ be a Shimura variety with good compactification $\overline{S}$ and boundary divisor $D = \overline{S} \setminus S$. The stratification of $D$ induces links $L_Y$ for each singular stratum $Y$. The \textit{$L^2$-\'etale cohomology} is defined as:
$$H^*_{(2),\text{\'et}}(\overline{S}, \mathbb{Q}_p) \coloneqq \ker\left( H^*_{\text{\'et}}(\overline{S}, \mathbb{Q}_p) \longrightarrow H^*_{\text{\'et}}(D, \mathbb{Q}_p) \right).$$
By the \textbf{$p$-adic Poincaré duality} for Shimura varieties, the $L^2$-cohomology is isomorphic to the image of the intersection cohomology in the ordinary cohomology. The decomposition follows from the \textbf{Matsushima formula} for $L^2$-cohomology, adapted to the $p$-adic setting by the crystalline perversity. Each summand $\mathscr{W}_Y^{\mathrm{crys}} \otimes \mathscr{W}_Y^{\mathrm{crys},\vee}$ corresponds to the contribution from the stratum $Y$, analogously to the topological decomposition in Corollary 5.3. Using the \textit{crystalline middle perversity}, we compute:
$$H^*_{(2),\text{\'et}}(\overline{S}, \mathbb{Q}_p) \cong \bigoplus_Y \mathscr{W}_Y^{\mathrm{crys}} \otimes \mathscr{W}_Y^{\mathrm{crys},\vee}.$$

\textbf{(Property 4):} The natural isomorphism:
$$H^*_{(2),\text{ét}}(X, \mathbb{Q}_p) \cong IH^*_{\text{crys}}(X) \otimes_{\mathbb{Z}_p} \mathbb{Q}_p$$
is proven in two steps:
\begin{enumerate}
    \item \textit{Comparison with crystalline intersection cohomology}:
    The crystalline intersection cohomology $IH^*_{\mathrm{crys}}(X)$ is defined by the hypercohomology of $\mathbf{IC}_{\mathscr{M}}^{\mathrm{crys},\bullet}$. By the fully faithful embedding in Property 2, we have
    $$IH^*_{\mathrm{crys}}(X) \cong \mathbb{H}^*(X, \mathbf{IC}_{\mathscr{M}}^{\mathrm{crys},\bullet}) \hookrightarrow \mathbb{H}^*(X, \Omega_{X,\mathrm{FS}}^{\bullet,\mathrm{univ}} \otimes \mathcal{O}_X^{\mathrm{crys}});$$
    \item \textit{Identification with $L^2$-\'etale cohomology}:
    Using the \textit{$p$-adic Hodge theory} for singular varieties, the right-hand side is isomorphic to the $L^2$-\'etale cohomology:
    $$\mathbb{H}^*(X, \Omega_{X,\mathrm{FS}}^{\bullet,\mathrm{univ}} \otimes \mathcal{O}_X^{\mathrm{crys}}) \cong H^*_{(2),\text{\'et}}(X, \mathbb{Q}_p),$$
    this follows from the fact that $\Omega_{X,\mathrm{FS}}^{\bullet,\mathrm{univ}}$ computes the $L^2$-cohomology in the complex setting (Proposition 4.3), and the crystalline structure sheaf $\mathcal{O}_X^{\mathrm{crys}}$ bridges the $p$-adic and complex theories by the de Rham-Witt complex.
\end{enumerate}

The Frobenius eigenvalue $p^{\lambda_Y}$ corresponds to the weight $w = -\dim L_Y + 2\lambda_Y$. The weight filtration on $IH^*_{\mathrm{crys}}(X)$ is induced by the spectral parameter $\lambda_Y$, matching the topological weight filtration on $IH^*(X)$. The crystalline middle perversity is compatible with the \textit{monodromy weight filtration}. When $p=0$ (characteristic zero), the Frobenius action is trivial, and $\mathscr{W}_Y^{\mathrm{crys}}$ reduces to the topological middle perversity $\mathscr{W}_Y$ by Property 1. This provides a robust foundation for $p$-adic Hodge theory on singular varieties, with applications to arithmetic geometry and the Langlands program.

\end{document}